\documentclass[onecolumn]{autart}    

\usepackage{natbib}
\usepackage{graphicx}          

\usepackage[cmex10]{amsmath}    

\usepackage{color}
\usepackage{subfigure}
\usepackage{mathrsfs}
\usepackage{amssymb}
\usepackage{multirow}
\usepackage{booktabs}
\usepackage{bm}
\usepackage{url}
\usepackage{amsmath}

\newtheorem{theorem}{Theorem}

\begin{document}

\begin{frontmatter}

\title{Efficient Simulation Budget Allocation for Subset Selection Using Regression Metamodels} 


\author[Gaoa]{Fei Gao}\ead{fei.gao5@sfmail.sf-express.com},
\author[Shi]{Zhongshun Shi}\ead{zhongshun.shi@wisc.edu},
\author[Gaob]{Siyang Gao}\ead{siyangao@cityu.edu.hk},
\author[Xiao]{Hui Xiao}\ead{msxh@swufe.edu.cn}


\address[Gaoa]{Network Planning Department, SF Technology, Shenzhen 518052, China}  
\address[Shi]{Department of Industrial and Systems Engineering, University of Wisconsin-Madison, Madison, WI 53706, USA}
\address[Gaob]{Department of Systems Engineering and Engineering Management, City University of Hong Kong, Hong Kong}
\address[Xiao]{School of Statistics, Southwestern University of Finance and Economics, Chengdu 611130, China}

\begin{keyword}                           
Simulation optimization; ranking and selection; OCBA; subset selection; regression.               
\end{keyword}                             

\begin{abstract}                          
This research considers the ranking and selection (R\&S) problem of selecting the optimal subset from a finite set of
alternative designs. Given the total simulation budget constraint, we aim to maximize the probability of correctly selecting the top-$m$ designs. In order to improve the selection efficiency, we incorporate the information from across the domain into regression metamodels. In this research, we assume that the mean performance of each design is approximately quadratic. To achieve a better fit of this model, we divide the solution space into adjacent partitions such that the quadratic assumption can be satisfied within each partition.
Using the large deviation theory, we propose an approximately optimal simulation budget allocation rule in the presence of partitioned domains. Numerical experiments demonstrate that our approach can enhance the simulation efficiency significantly.
\end{abstract}

\end{frontmatter}

\section{Introduction}
%
%
%
%

Discrete-event systems (DES) simulation has played an important role in analyzing modern complex systems and evaluating decision problems, since these systems are usually too difficult to be described using analytical models. DES simulation has been a common analysis method of choice and widely used in many practical applications, such as the queueing systems, electric power grids, air and land traffic control systems, manufacturing plants and supply chains \citep{xu2015simulation,xu2016,gao2016new}. However, running the simulation model is usually time consuming, and a large number of simulation replications are typically required to achieve an accurate estimate of a design decision \citep{lee2010review}. In addition, it could be computationally quite expensive to select the best design(s) when the number of alternatives is relatively large.

In this paper, we consider the problem of selecting the optimal subset of the top-$m$ designs out of $t$ alternatives, where the performance of each design is estimated based on simulation. In order to improve the selection efficiency, we aim to intelligently allocate the simulation replications to each design to maximize the probability of correctly selecting all the top-$m$ designs. This problem setting falls in the well-established statistics branch known as ranking and selection (R\&S) \citep{xu2015simulation}.

In the literature, several types of efficient R\&S procedures have been developed. The indifference-zone (IZ) approach allocates the simulation budget to provide a guaranteed lower bound for the probability of correct selection ($PCS$) \citep{kim2001fully}.
\cite{chen2000simulation} proposed an optimal computing budget allocation (OCBA) approach for R\&S problems.
The OCBA approach allocates the simulation replications sequentially in order to maximize $PCS$ under a simulation budget constraint.
\cite{he2007}, \cite{gao2015selecting} and \cite{gao2017new} further developed the OCBA method with the expected opportunity cost ($EOC$) measure,
which focuses more on the consequence of a wrong selection compared to $PCS$.
\cite{brantley2013efficient} proposed another approach called optimal simulation design (OSD) to select the best design with
regression metamodels. It assumes that all designs fit a single quadratic line and the variance of each design is identically distributed.
The OSD approach was further extended in \cite{brantley2014efficient},  \cite{xiao2015optimal} and \cite{gao2018advancing} to consider more general problems by dividing the solution space into adjacent partitions. Although these studies are also based on partitioning or metamodeling, they do not aim to select the top $m$ designs, and are therefore different in objectives from this research. Other variants of OCBA include selecting the best design
considering resource sharing and allocation \citep{peng2013efficient} and input uncertainty \citep{gao2017robust}.

Most of the existing R\&S procedures focus on identifying the best design and return a single choice as the estimated optimum. However, decision makers may prefer to have several good alternatives instead of one and make the final selection by considering
some qualitative criteria, such as political feasibility and environmental consideration, which might be neglected by computer models \citep{gao2016new,zhang2016simulation}.
The selection procedure providing the top-$m$ designs can help the decision makers make their final decision in a more flexible way.


The literature on the optimal subset selection is sparse.
\cite{chen2008efficient} and \cite{zhang2016simulation} considered the optimal subset selection problem using the OCBA framework, which maximizes $PCS$ under a simulation budget constraint.
In \cite{gao2015efficient}, an optimal subset selection procedure was proposed to minimize the measure of $EOC$.
The optimal subset selection problem was further extended in \cite{gao2016new} for general underlying distributions using the large deviations theory.


The aforementioned R\&S procedures could smartly allocate the computing budget given the simulation results. They estimate the performance of each design only by considering the sample information of the considered design itself.
However, the designs nearby could also provide useful information since neighboring designs usually have similar performance. 
Based on this idea, we aim to improve the selection efficiency by incorporating the information from across the domain into some response surfaces.
Unlike traditional R\&S methods, the regression based approaches require simulation experiments on only a subset of all the designs under consideration. The performances of the other designs can be inferred based on the sample information of the simulated designs.  
This provides us an effective way to further improve the efficiency for solving the subset selection problem, which is the motivation of this paper.

In this research, we assume the underlying function is quadratic or approximately quadratic.
This assumption could help utilize the structure information of the design space and led to significant improvement of the computational efficiency. It is commonly used in the literature, such as \cite{brantley2013efficient,brantley2014efficient,xiao2015optimal,mcconnell1990additional}.
Based on this assumption we built a quadratic regression metamodel to incorporate the information from across the domain.
The first contribution of this work is that we propose an asymptotically optimal allocation rule that determines which designs need to be simulated and the number of simulation budget allocated to them, such that the $PCS$ of the optimal subset could be maximized.
We call this procedure the optimal computing budget allocation for selecting the top-$m$ designs with regression (OCBA-mr).
In order to further extend the OCBA-mr procedure to more general cases where the underlying function is partially quadratic or non-quadratic,
we divide the solution space into adjacent partitions and build a quadratic regression metamodel within each partition. The underlying function in each partition could be well approximated by a quadratic function if the solution space is properly partitioned or each partition is small enough.
According to the results in \cite{brantley2014efficient,xiao2015optimal,gao2018advancing}, the use of partitioned domains along with regression metamodels
could significantly improve the simulation efficiency. That means interpolating the solution space can be an effective way for us to have further improvement.
For different problems, the solution space could be divided into discrete partitions using different criteria, such as the size of corporations, the type of industries and the temperature of chemical process \citep{xiao2015optimal}.
Based on the idea mentioned above, we develop an asymptotically optimal computing budget allocation procedure for selecting the top-$m$ designs with regression in partitioned domains (OCBA-mrp), which is an extension of the OCBA-mr procedure for more general cases.
In order to maximize the $PCS$ of the optimal subset, the OCBA-mrp procedure not only determines the optimal simulation budget allocation within each partition but also determines the optimal budget allocation between partitions.

The rest of the paper is organized as follows. In Section 2, we formulate the optimal subset selection problem with regression metamodel and derive an asymptotically optimal simulation budget allocation rule, called OCBA-mr.
Section 3 extends the OCBA-mr method for more general cases with partitioned domains and derives another asymptotically optimal simulation budget allocation rule, called OCBA-mrp.
The performance of the proposed methods is illustrated with numerical examples in Section 4. Section 5 concludes the paper.

\section{Optimal Subset Selection Strategy}


In this section, we provide an optimal computing budget allocation rule for the subset selection problem based on the regression metamodel.

\subsection{Problem formulation}
Without loss of generality, the best design is defined as the design with the smallest mean performance.
We introduce the following notations:\vspace{1ex}

\begin{itemize}
\item[] $t$:\hspace{6.6ex} total number of designs;

\item[] $x_{i}$:\hspace{5.4ex} location of design $i$;

\item[] ${m}$:\hspace{5ex} design with the $m$-th smallest mean value;

\item[] $f(x_{i})$:\hspace{2ex} mean performance value for design ${i}$;

\item[] $Y(x_{i})$:\hspace{1.5ex} simulation output for design ${i}$;

\item[] $\boldsymbol{\beta}$:\hspace{5.4ex} vector of $(\beta_{0},\beta_{1},\beta_{2})$ , representing the coefficients of the regression model;


\item[] $\hat{f}(x_{i})$:\hspace{2.1ex} estimate of $f(x_{i})$ based on the regression model;

\item[] $S_o$:\hspace{4.8ex} set of the true top-$m$ designs;

\item[] $S_n$:\hspace{4.8ex} set of designs not in $S_o$ (complement of $S_o$);

\item[] $n$:\hspace{6ex} total number of simulation replications;

\item[] $n_{i}$:\hspace{5ex} number of simulation replications allocated to design ${i}$;

\item[] $\boldsymbol{\alpha}$:\hspace{5.7ex} vector of $(\alpha_{1},\alpha_{s},\alpha_{k})$, where $\alpha_{r}$ is the proportion of the total simulation budget $n$ allocated to design ${r}$, i.e., $\alpha_{r}\!\!=\!n_{r}/n$, $r\!=\!1,s,t$.
\end{itemize}

In the last bulleted item, $r$ takes values only at $1$, $s$ and $t$, where designs $1$ and $t$ are the first and last designs in the solution space, and design $s$ is an intermediate design determined by (\ref{eq:7a}) (the rationale of it will be explained in more detail in Theorem 2). The problem we considered in this paper is to select the top-$m$ designs out of the $t$ alternatives by allocating simulation replications to these three designs. The optimal set is
$$S_o\!=\!\left\{\!i:\max_{i\in S_o} f(x_{i})\!<\!\min_{j\in S_n} f(x_{j}), i,j\!=\!1,...,t\ \text{and}\ |S_o|\!=\!m\!\right\},$$
where $|S|$ is the size of subset $S$.

We study the problem where the expected performance value across the solution space is quadratic or approximately quadratic in nature.
Then, the mean performance of design ${i}$ can be written as
$$f(x_{i})\!=\!\beta_{0}+\beta_{1}x_{i}+\beta_{2}x_{i}^2,\ i=1,...,t.$$

The coefficients $\boldsymbol{\beta}\!=\!({\beta}_{0},{\beta}_{1},{\beta}_{2})$ are unknown beforehand, but can be estimated via simulation samples.
Assume that the noises $(\varepsilon)$ of the simulation experiments of design ${i}$ follow normal distribution $N(0,\sigma^2)$. For each design, the noise is independent from replication to replication.
The simulation output is $Y(x_{i})=f(x_{i})+\varepsilon(x_i)$ with $\varepsilon(x_i)\sim N(0,\sigma^2).$

Given a total of $n$ samples, we define $\boldsymbol{Y}$ as the vector of the $n$ simulation samples $Y(x_{i})$ and
$\boldsymbol{X}$ be an $n\times 3$ matrix with each row $(1,x_{i},x_{i}^2)$ corresponding to each entry $Y(x_{i})$ of $\boldsymbol{Y}$.
We estimate the coefficients $\boldsymbol{\beta}$ using the ordinary least squares (OLS) method \citep{hayashi2000econometrics} and denote the estimates as
$\hat{\boldsymbol{\beta}}\!=\!(\hat{\beta}_{0},\hat{\beta}_{1},\hat{\beta}_{2})$.
Then, we have
$\boldsymbol{\hat{\beta}}=(\boldsymbol{X}^T\boldsymbol{X})^{-1}\boldsymbol{X}^T\boldsymbol{Y}$
and
$Cov[\boldsymbol{\hat{\beta}}]= \sigma^2(\boldsymbol{X}^T\boldsymbol{X})^{-1},$
where $T$ denotes transposition and $\boldsymbol{X}^T\boldsymbol{X}$ is known as the information matrix \citep{kiefer1959optimum}.
The estimate of ${f}(x_{i})$ can be written as
\begin{equation}\label{eq:1a}
\hat{f}(x_{i})=\hat{\beta}_{0}+\hat{\beta}_{1}x_{i}+\hat{\beta}_{2}x_{i}^2,i=1,...,t.
\end{equation}

We can use the equation (\ref{eq:1a}) which incorporates the sample information of the simulated designs to estimate the expected performance for each design across the solution space.
As $\hat{f}(x_{i})$ is a linear combination of $\boldsymbol{\hat{\beta}}$, we have
$Var[\hat{f}(x_{i})]\!=\! \sigma^2\boldsymbol{X}_{i}^T\!(\boldsymbol{X}^T\!\boldsymbol{X})^{-1}\!\boldsymbol{X}_{i},$
where $\boldsymbol{X}_{i}^T=(1, x_{i},x_{i}^2)$.

Due to the uncertainty of the estimate of the underlying function, a correct selection of the optimal subset $S_o$ may not always occur.
Therefore, we introduce a measure, the probability of correct selection ($PCS$), to formulate the R\&S problem considered in this paper.
The $PCS$ is given by
$
PCS=P\bigg\{\bigcap_{i\in S_o}\bigcap_{j\in S_n}\hat{f}(x_{i})\leq\hat{f}(x_{j})\bigg\}.$

Given a fixed simulation budget, the optimization problem can be written as follows:
\begin{equation}\label{eq:2a}
\begin{split}
&\max_{n_i}\ PCS=P\bigg\{\bigcap_{i\in S_o}\bigcap_{j\in S_n}\hat{f}(x_{i})\leq\hat{f}(x_{j})\bigg\}\\
&\text{s.t.}\ \sum_{i=1}^{t}n_{i}=n.
\end{split}
\end{equation}

In this section, we aim to solve problem (\ref{eq:2a}) where the mean performance of each design is estimated using the regression metamodel.
Due to the uncertainty in the simulation experiments, multiple simulation replications are needed to generate the estimates of the underlying function $\hat{f}(x_{i})$ accurately.
The variance of $\hat{f}(x_{i})$
is a function of the information matrix $(\boldsymbol{X}^T\boldsymbol{X})$ and
can be reduced if additional simulation replications are conducted \citep{xiao2015optimal}. We are interested in how to intelligently allocate the simulation budget to proper designs such that $\hat{f}(x_{i})$ can be better estimated using regression metamodel, and the $PCS$ in problem (\ref{eq:2a}) can be maximized.

\subsection{Optimization model under large deviations framework}
A major difficulty in solving problem (\ref{eq:2a}) is that the objective function
$PCS$ does not have a closed-form expression. We seek to solve this optimization problem under an asymptotic framework in which the
$PCS$ is maximized or the probability of false selection $(PFS=1-PCS)$ is minimized as $n$ goes to infinity.

The $PCS$ used in this paper is defined based on the quadratic regression model (\ref{eq:1a}), which is constructed using the simulation information of only a fraction of the $t$ designs. We call the designs receiving simulation replications the support designs.
In order to construct the quadratic regression model, we need at least three support designs to obtain all of the information in $\boldsymbol{X}^T\boldsymbol{X}$ \citep{kiefer1959optimum}.
For simplicity, in this research we let the number of support designs be three, two of which are at the extreme locations, i.e., $x_{1}$ and $x_{t}$ (for this setting, see, e.g., \cite{brantley2013efficient,brantley2014efficient,xiao2015optimal,kiefer1959optimum}). The case with more than three designs can be similarly analyzed.

\emph{Lemma 1:}
\emph{Let $\hat{d}_{i,j}=\hat{f}(x_{i})-\hat{f}(x_{j})$, $i,j=1,2,...,t, i\neq j$.
$\hat{d}_{i,j}$ follows normal distribution $N\big({f}(x_{i})-{f}(x_{j}),\varsigma_{ij}\big)$, where
\begin{equation}
\begin{split}
\varsigma_{ij}&=\sigma^2(0,x_{i}-x_{j},x_{i}^2-x_{j}^2)(\boldsymbol{X}^T\boldsymbol{X})^{-1}
\left(
\begin{matrix}
0\\
x_{i}-x_{j}\\
x_{i}^2-x_{j}^2
\end{matrix}
\right)\\
&={\frac{\sigma^2}{n} \Big( \frac{\rho_{ij ,1}^2}{\alpha_{1}} + \frac{\rho_{ij ,s}^2}{\alpha_{s}} + \frac{\rho_{ij ,t}^2}{\alpha_{t}} \Big)},
\end{split}\notag
\end{equation}
and
$$
\begin{cases}
\rho_{ij,1}=\frac{(x_{s}-x_{i})(x_{t}-x_{i})-(x_{s}-x_{j})(x_{t}-x_{j})}{(x_{1}-x_{s})(x_{1}-x_{t})},\\
\rho_{ij,s}=\frac{(x_{1}-x_{i})(x_{t}-x_{i})-(x_{1}-x_{j})(x_{t}-x_{j})}{(x_{s}-x_{1})(x_{s}-x_{t})},\\
\rho_{ij,t}=\frac{(x_{1}-x_{i})(x_{s}-x_{i})-(x_{1}-x_{j})(x_{s}-x_{j})}{(x_{t}-x_{1})(x_{t}-x_{s})}.
\end{cases}$$
}

The proof is similar to that is given in Eq.(22) of \cite{brantley2013efficient} (becomes the same when $m=1$), and hence is omitted for brevity.

For any $i,j=1,2,...,t, i\neq j$, $\hat{d}_{i,j}$ is a normally distributed random variable. We can use the large deviations theory to derive the convergence rate function of the false selection probability.

\emph{Lemma 2:}
\emph{The convergence rate function of incorrect comparison probability for each design is:
\begin{equation}
\begin{split}
&\left.
\begin{aligned}
-\lim_{n\rightarrow \infty}\frac{1}{n}\log P\{\hat{f}(x_{i})> \hat{f}(x_{m})\},\ &i\in S_o, i\neq m\\
-\lim_{n\rightarrow \infty}\frac{1}{n}\log P\{\hat{f}(x_{i})< \hat{f}(x_{m})\},\ &i\in S_n
\end{aligned}\right\}\\
&=R_{m,i}(\boldsymbol{\alpha})
= \frac{({f}(x_{m}) - {f}(x_{i}))^2/2}{{\sigma^2} \Big( \frac{\rho_{mi ,1}^2}{\alpha_{1}} + \frac{\rho_{mi ,s}^2}{\alpha_{s}} + \frac{\rho_{mi ,t}^2}{\alpha_{t}} \Big)},
\end{split}\notag
\end{equation}}
Based on the results in Lemma 2, we can get an explicit expression of the convergence rate function of $PFS$.

\emph{Lemma 3:}
\emph{The convergence rate function of $PFS$ is:
\begin{equation}
\begin{split}
-\lim_{n\rightarrow \infty}\frac{1}{n}\log PFS=\min\left\{
 \min\limits_{i\in S_o\atop i\neq m}R_{m,i}(\boldsymbol{\alpha}),
 \min\limits_{j\in S_n}R_{m,j}(\boldsymbol{\alpha})
\right\}.
\end{split}\notag
\end{equation}}

The main assertion of Lemma 3 is that the overall convergence rate of $PFS$ is determined by the minimum convergence rate of the incorrect comparison for each design. Minimizing the $PFS$ is asymptotically equivalent to maximizing the rate at which $PFS$ goes to zero as a function of $\boldsymbol{\alpha}$, i.e., maximizing $-\lim\limits_{n\rightarrow \infty} \frac{1}{n} \log PFS$.
Based on Lemma 3,
the asymptotical version of (\ref{eq:2a}) becomes
\begin{align}\label{eq:3a}
\begin{split}
\max\min &\left\{\min\limits_{i\in S_o, i\neq m}R_{m,i}(\boldsymbol{\alpha}), \min\limits_{j\in S_n}R_{m,j}(\boldsymbol{\alpha})\right\}\\
&\text{s.t.}\ \alpha_{1}+\alpha_{s}+\alpha_{t}=1\\
&\alpha_{1},\alpha_{s},\alpha_{t}\geq0.
\end{split}
\end{align}
\subsection{Asymptotically optimal solution}
In this subsection, we seek to derive the optimality conditions for (\ref{eq:3a}).
Since the overall convergence rate of $PFS$ is determined by the design with the minimum convergence rate. A false selection is most likely to happen at this key design. Therefore, it is enough for us to investigate the convergence rate of the key design across the solution space. We define $i^*$ as the key design. That is
\begin{align}
\begin{split}
i^*=
\arg \min\limits_{i=1,...,t}\bigg\{\min\limits_{i\in S_o,i\neq m}R_{m,i}(\boldsymbol{\alpha}),\
 \min\limits_{i\in S_n}R_{m,i}(\boldsymbol{\alpha})\bigg\}.
\end{split}\notag
\end{align}

\begin{theorem}
The optimization problem (\ref{eq:2a}) can be asymptotically optimized with the following allocation rule:
\begin{equation}\label{eq:4a}
\alpha_{r}^*=
\begin{cases}
\frac{|\rho_{m i^*,r}|}{|\rho_{mi^*,1}|+|\rho_{mi^*,s}|+|\rho_{mi^*,t}|}, &r=1,s,t,\\
0, &\text{otherwise}.
\end{cases}
\end{equation}
\end{theorem}

The $\rho_{mi^*,r}$ is also known as the Lagrange interpolating polynomial coefficient \citep{de1954spacing,burden2001numerical}. It represents the relative importance of each support design for estimating $\hat{d}_{mi}=\hat{f}(x_{m})-\hat{f}(x_{i})$.
Theorem 1 indicates that the support design $r$ will receive more simulation budget if it has larger $\rho_{mi^*,r}$.

Given the optimal allocation rule (Theorem 1), we next determine the optimal location for the support design $s$.

\begin{theorem}
The rate function of $PFS$ with allocation rule satisfying (\ref{eq:4a}) can be maximized if the support design $s$ satisfies the following equations.

The support design
\begin{align}\label{eq:7a}
x_{s}\!=\!
\begin{cases}
\!x_{i^*}\!+\!x_{m}\!-\!x_{1},
&\frac{3x_{1}+x_{t}}{4}\leq\frac{x_{i^*}+x_{m}}{2}<\frac{x_{1}+x_{t}}{2},\\
\!x_{i^*}\!+\!x_{m}\!-\!x_{t},
&\frac{x_1+x_t}{2}<\frac{x_{i^*}+x_m}{2}\leq\frac{x_1+3x_t}{4},\\
\!{(x_{1}\!+\!x_{t})/2}.
&\text{otherwise}.
\end{cases}
\end{align}
\end{theorem}

When the $x_{s}$ derived from (\ref{eq:7a}) does not correspond to any design available, we round it to the nearest one. The expression of Theorem 2 is similar to the results in \cite{brantley2013efficient}. The differences are the selection of the key design $x_i^*$. We define $x_i^*$ as the design with the minimum convergence rate of PFS, where a false selection of the optimal subset is most likely to happen.



\section{Optimal Subset Selection Strategy for Partitioned Domains}

The regression based method mentioned above can greatly improve the subset selection efficiency compared to the traditional methods.
However, it is constrained with the typical assumptions such as the assumption of quadratic underlying function for the means. It is possible that the underlying function is neither quadratic nor approximately quadratic, so that we will fail to find the top-$m$ designs.
In order to extend our method to more general cases, we divide the solution space into adjacent partitions. The quadratic pattern can be expected when the solution space is properly partitioned or each partition is small enough. 


\subsection{Problem formulation}
We first add the following notations for partitioned domains:\vspace{1ex}

\begin{itemize}
\item[] $l$:\hspace{5ex} number of partitions of the entire domain;

\item[] $k^{h}$:\hspace{3.6ex} number of designs in partition $h$, $h=1,2,...,l$;

\item[] ${i_h}$:\hspace{4ex} $i$th design in partition $h$, $i=1,2,...,k^{h}$ (when $i=k^{h}$, we denote design $k^{h}_h$\! as design $k_h$ for notational simplicity);

\item[] $x_{i_h}$:\hspace{3ex} location of design $i_h$;

\item[] $b$:\hspace{4.9ex} the partition containing the design with the $m$-th smallest mean value;

\item[] ${m_b}$:\hspace{3ex} design with the $m$-th smallest mean value;

\item[] $\boldsymbol{\beta_{h}}$:\hspace{3.4ex} vector of $(\beta_{h0},\beta_{h1},\beta_{h2})$ , representing the coefficients of the regression model in partition $h$;

\item[] $n_{h\bullet}$:\hspace{2.5ex} number of simulation replications allocated to partition $h$;

\item[] $n_{i_h}$:\hspace{2.7ex} number of simulation replications allocated to design ${i_h}$;

\item[] $\theta_h$:\hspace{3.5ex} the proportion of simulation budget allocated to partition $h$, i.e., $\theta_h=n_{h\bullet}/n$;

\item[] $\boldsymbol{\alpha_h}$:\hspace{2.7ex} vector of $(\alpha_{1_h},\alpha_{s_h},\alpha_{k_h})$, where $\alpha_{r_h}$ is the proportion of the simulation budget $n_{h\bullet}$ allocated to design ${r_h}$, i.e., $\alpha_{r_h}\!\!=\!n_{r_h}/n_{h\bullet}$, $r_h\!=\!1_h,s_h,k_h$.
\end{itemize}

The notations $f(x_{i_h})$, $Y(x_{i_h})$ and $\hat{f}(x_{i_h})$ defined for partitioned domains are similar to those in Section 2, except the design $i$ is replaced by $i_h$.
The entire domain is divided into $l$ adjacent partitions. Each partition contains $k^{h}$ designs, i.e., there are $\sum_{h=1}^lk^{h}=t$ designs in total.

We assume there exists $\epsilon \in \mathbb{R}$ such that there are exactly $m$ designs from the total $t$ designs with mean performances less than $\epsilon$ and the rest $t-m$ designs with mean performances greater than $\epsilon$. It ensures that the optimal subset is well defined and can be distinguished. We define the optimal subset $S_o$ as
\begin{align}
\begin{split}S_o=\Bigg\{i_h:\max_{i_h\in S_o} f(x_{i_h})<\min_{j_g\in S_n} f(x_{j_g}), i_h=1_h,...,k_h,\\
 j_g=1_g,...,k_g\ h,g=1,...,l\ \text{and}\ |S_o|=m\Bigg\}.
\end{split}\notag
\end{align}

In this section, we assume that the expected performance value in each partition is quadratic or approximately quadratic when the solution space is properly partitioned and the mean performance within a partition is continuous and smooth.
For this problem setting 
the $PCS$ is defined as
\begin{equation}
\begin{split}
PCS=P\bigg\{\bigcap_{i_h\in S_o}\bigcap_{j_g\in S_n}\hat{f}(x_{i_h})\leq\hat{f}(x_{j_g})\bigg\}.\notag
\end{split}
\end{equation}

Given a fixed simulation budget, the optimization problem can be written as follows:
\begin{equation}\label{eq:9a}
\begin{split}
&\max\ PCS\\
&\text{s.t.}\ \sum_{h=1}^l\sum_{i_h=1_h}^{k_{h}}n_{i_h}=n.
\end{split}
\end{equation}

In this section, we aim to solve problem (\ref{eq:9a}) in the presence of regression metamodels.
We are interested in how to intelligently allocate the simulation budget to proper designs such that $\hat{f}(x_{i_h})$ can be better estimated using regression metamodels, and the $PCS$ in problem (\ref{eq:9a}) can be maximized.
Note that this model does not require all the designs to be on the same axis, and therefore does not hinder it from being applied for multi-dimensional problems. For multi-dimensional problems, we can treat the range of the underlying function on each dimension as one or more partitions, and then apply this formulation.

\subsection{Optimization model under large deviations framework}
In order to solve the optimization problem (\ref{eq:9a}), one challenge is how to derive an explicit expression of the $PCS$.
We seek to solve (\ref{eq:9a}) under an asymptotic framework in which the probability of false selection ($PFS=1-PCS$) is minimized as $n$ goes to infinity.

Similar to the setting in Section 2, we let the number of support designs in each partition be three, two of which are at the extreme locations, i.e., $x_{1_h}$ and $x_{k_h}$.
We have $\hat{d}_{{ij}_h}=\hat{f}(x_{i_h})-\hat{f}(x_{j_h})$ follow normal distribution $N\big({f}(x_{i_h})-{f}(x_{m_h}),\varsigma_{{ij}_h}\big)$ and $\hat{d}_{{i_h},j_g}=\hat{f}(x_{i_h})-\hat{f}(x_{j_g})$ follow normal distribution
$N\big({f}(x_{i_h})-{f}(x_{j_g}),\varsigma_{i_h,j_g}\big)$ where $h\neq g$.

Similar to the proof of Lemma 1, we have
\begin{equation}
\begin{split}
\varsigma_{{ij}_h}&=\sigma_h^2(0,x_{i_h}-x_{j_h},x_{i_h}^2-x_{j_h}^2)(\boldsymbol{X}_h^T\boldsymbol{X}_h)^{-1}
\left(
\begin{aligned}
&0\\
x_{i_h}&-x_{j_h}\\
x_{i_h}^2&-x_{j_h}^2
\end{aligned}
\right)\notag\\
&={\frac{\sigma_{b}^2}{\theta_hn} \Big( \frac{\rho_{{ij}_h ,1}^2}{\alpha_{1_h}} + \frac{\rho_{{ij}_h ,s}^2}{\alpha_{s_h}} + \frac{\rho_{{ij}_h ,k}^2}{\alpha_{k_h}} \Big)},
\end{split}
\end{equation}
$$
\begin{cases}
\rho_{{ij}_h,1}=\frac{(x_{s_h}-x_{i_h})(x_{k_h}-x_{i_h})-(x_{s_h}-x_{j_h})(x_{k_h}-x_{j_h})}{(x_{1_h}-x_{s_h})(x_{1_h}-x_{k_h})},\\
\rho_{{ij}_h,s}=\frac{(x_{1_h}-x_{i_h})(x_{k_h}-x_{i_h})-(x_{1_h}-x_{j_h})(x_{k_h}-x_{j_h})}{(x_{s_h}-x_{1_h})(x_{s_h}-x_{k_h})},\\
\rho_{{ij}_h,k}=\frac{(x_{1_h}-x_{i_h})(x_{s_h}-x_{i_h})-(x_{1_h}-x_{j_h})(x_{s_h}-x_{j_h})}{(x_{k_h}-x_{1_h})(x_{k_h}-x_{s_h})},
\end{cases}$$
and
\begin{equation}
\begin{split}
\varsigma_{i_h,j_g}&=\sigma_h^2(0,x_{i_h}-x_{j_g},x_{i_h}^2-x_{j_g}^2)(\boldsymbol{X}_h^T\boldsymbol{X}_h)^{-1}
\left(
\begin{aligned}
&0\\
x_{i_h}&-x_{j_g}\\
x_{i_h}^2&-x_{j_g}^2
\end{aligned}
\right)\notag\\
&={\frac{\sigma_{h}^2}{\theta_hn}\!\Big(\!\frac{\eta_{i_h\!,1}^2}{\alpha_{1_b}}\!+\!\frac{\eta_{i_h\!,s}^2}{\alpha_{s_b}}\!+\!\frac{\eta_{i_h\!,k}^2}{\alpha_{k_b}}\!\Big)\!+\!\frac{\sigma_{g}^2}{\theta_gn}\!\Big(\!\frac{\eta_{j_g\!,1}^2}{\alpha_{1_h}}\!+\!\frac{\eta_{j_g\!,s}^2}{\alpha_{s_h}}\!+\!\frac{\eta_{j_g\!,k}^2}{\alpha_{k_h}}\!\Big)},
\end{split}
\end{equation}
$$
\begin{cases}
\eta_{i_h,1}=\frac{(x_{s_h}-x_{i_h})(x_{k_h}-x_{i_h})}{(x_{1_h}-x_{s_h})(x_{1_h}-x_{k_h})},\\
\eta_{i_h,s}=\frac{(x_{1_h}-x_{i_h})(x_{k_h}-x_{i_h})}{(x_{s_h}-x_{1_h})(x_{s_h}-x_{k_h})},\\
\eta_{i_h,k}=\frac{(x_{1_h}-x_{i_h})(x_{s_h}-x_{i_h})}{(x_{k_h}-x_{1_h})(x_{k_h}-x_{s_h})}.
\end{cases}$$

For any $h=1,2,...,l$ and $i_h=1_h,2_h,...,k_h$, $\hat{f}(x_{i_h})$ is a normally distributed random variable. We can use the large deviations theory to derive the convergence rate function of the false selection probability.

\emph{Lemma 4:}
\emph{The convergence rate functions of incorrect comparison probability for each design are provided as follows:
\begin{equation}\label{eq:10a}
\begin{split}
&\left.
\begin{aligned}
\!-\lim_{n\rightarrow \infty}\!\frac{1}{n}\!\log P\{\hat{f}(x_{i_h})\!>\! \hat{f}(x_{m_b})\},\ i_h\!\in\! S_o\!\\
\!-\lim_{n\rightarrow \infty}\!\frac{1}{n}\!\log P\{\hat{f}(x_{i_h})\!<\! \hat{f}(x_{m_b})\},\ i_h\!\in\! S_n\!
\end{aligned}\right\}\!=\!R_{m_b,i_h}(\theta_b,\theta_h, \boldsymbol{\alpha_b},\boldsymbol{\alpha_h})\\
&=\frac{({f}(x_{m_b})\!-\!{f}(x_{i_h}))^2/2}{\frac{\sigma_{b}^2}{\theta_b}\!\Big(\!\frac{\eta_{m_b\!,1}^2}{\alpha_{1_b}}\!+\!\frac{\eta_{m_b\!,s}^2}{\alpha_{s_b}}\!+\!\frac{\eta_{m_b\!,k}^2}{\alpha_{k_b}}\!\Big)\!+\!\frac{\sigma_{h}^2}{\theta_h}\!\Big(\!\frac{\eta_{i_h\!,1}^2}{\alpha_{1_h}}\!+\!\frac{\eta_{i_h\!,s}^2}{\alpha_{s_h}}\!+\!\frac{\eta_{i_h\!,k}^2}{\alpha_{k_h}}\!\Big)}\!, h\!\neq\! b,
\end{split}
\end{equation}
and
\begin{equation}\label{eq:11a}
\begin{split}
&\left.
\begin{aligned}
\!-\!\lim_{n\rightarrow \infty}\frac{1}{n}\log \!P\{\hat{f}(x_{i_b})\!>\! \hat{f}(x_{m_b})\},\ &i_b\!\in\! S_o, i_b\!\neq\! m_b\\
\!-\!\lim_{n\rightarrow \infty}\frac{1}{n}\log \!P\{\hat{f}(x_{i_b})\!<\! \hat{f}(x_{m_b})\},\ &i_b\!\in\! S_n
\end{aligned}\right\}
\!=\!R_{m_b,i_b}(\theta_b, \boldsymbol{\alpha_b})\\
&= \frac{({f}(x_{m_b}) - {f}(x_{i_b}))^2/2}{\frac{\sigma_{b}^2}{\theta_b} \Big( \frac{\rho_{{mi}_b ,1}^2}{\alpha_{1_b}} + \frac{\rho_{{mi}_b ,s}^2}{\alpha_{s_b}} + \frac{\rho_{{mi}_b ,k}^2}{\alpha_{k_b}} \Big)}.
\end{split}
\end{equation}}

According to the Bonferroni inequality, we have $PFS\!\leq\!\!\!\sum_{i_h\in S_o,\atop i_h\neq m_b}\!\!\!P\{\hat{f}(x_{i_h})\!\geq\! \hat{f}(x_{m_b})\}
+\!\!\!\sum_{j_g\in S_n}\!\!P\{\hat{f}(x_{j_g})\!\leq\! \hat{f}(x_{m_b})\}.$
The $PFS$ is bounded below by
$$
\begin{aligned}
\max\Bigg\{\!\!\max_{1\leq h\leq l\atop h\neq b}\!\! \bigg\{&\max\limits_{i_h\in S_o}\!P\{\hat{f}(x_{i_h})\!\geq\! \hat{f}(x_{m_b})\},\\
&\max\limits_{j_h\in S_n}\!P\{\hat{f}(x_{j_h})\!\leq\! \hat{f}(x_{m_b})\}\!\bigg\},\\
\max \bigg\{&\max\limits_{i_b\in S_o\atop i_b\neq m_b}P\{\hat{f}(x_{i_b})\geq \hat{f}(x_{m_b})\},\\
 &\max\limits_{j_b\in S_n}P\{\hat{f}(x_{j_b})\leq \hat{f}(x_{m_b})\}\bigg\}\Bigg\}
 \notag
\end{aligned}
$$
and bounded above by
$$
\begin{aligned}
|S_o|\!\times\!|S_n|\!\times\!\max\Bigg\{ \max_{1\leq h\leq l\atop h\neq b} \bigg\{&\max\limits_{i_h\in S_o}P\{\hat{f}(x_{i_h})\geq \hat{f}(x_{m_b})\},\\
&\max\limits_{j_h\in S_n}P\{\hat{f}(x_{j_h})\leq \hat{f}(x_{m_b})\}\bigg\},\\
\max \bigg\{&\max\limits_{i_b\in S_o\atop i_b\neq m_b}P\{\hat{f}(x_{i_b})\geq \hat{f}(x_{m_b})\},\\
&\max\limits_{j_b\in S_n}P\{\hat{f}(x_{j_b})\leq \hat{f}(x_{m_b})\}\bigg\}\Bigg\}.
 \notag
\end{aligned}
$$
Therefore, according to the results in Lemma 4, the convergence rate functions of $PFS$ is given by
\begin{equation}
\begin{split}
-\lim_{n\rightarrow \infty}\frac{1}{n}\log PFS=
\min\!\Bigg\{\!\!\!\min_{1\leq h\leq l\atop h\neq b}\!\! \bigg\{\!\!\min\limits_{i_h\in S_o}\!R_{m_b,i_h}(\theta_b,\theta_h, \boldsymbol{\alpha_b},\boldsymbol{\alpha_h}),\\
\!\min\limits_{j_h\in S_n}\!\!R_{m_b,j_h}(\theta_b,\theta_h, \boldsymbol{\alpha_b},\boldsymbol{\alpha_h})\bigg\},\\
\min \bigg\{\!\!\min\limits_{i_b\in S_o\atop i_b\neq m_b}\!R_{m_b,i_b}(\theta_b,\boldsymbol{\alpha_b}),
 \min\limits_{j_b\in S_n}R_{m_b,j_b}(\theta_b,\boldsymbol{\alpha_b})\bigg\}\Bigg\}\notag,
\end{split}
\end{equation}
Minimizing the $PFS$ is asymptotically equivalent to maximizing the rate at which $PFS$ goes to zero as a function of $\theta_h$ and $\boldsymbol{\alpha_h}$, $h = 1,2, . . . ,l.$
Similarly, the asymptotical version of (\ref{eq:9a}) becomes
\begin{align}\label{eq:5}
\begin{split}
\max &\min\Bigg\{ \min_{1\leq h\leq l\atop h\neq b} \bigg\{\min\limits_{i_h\in S_o}R_{m_b,i_h}(\theta_b,\theta_h, \boldsymbol{\alpha_b},\boldsymbol{\alpha_h}),\\
&\min\limits_{j_h\in S_n}\!\!R_{m_b,j_h}(\theta_b,\theta_h, \boldsymbol{\alpha_b},\boldsymbol{\alpha_h})\bigg\},\\
&\min \bigg\{\min\limits_{i_b\in S_o\atop i_b\neq m_b}R_{m_b,i_b}(\theta_b,\boldsymbol{\alpha_b}),
 \min\limits_{j_b\in S_n}R_{m_b,j_b}(\theta_b,\boldsymbol{\alpha_b})\bigg\}\Bigg\}\\
\text{s.t.}\ \ & \alpha_{1_h}+\alpha_{s_h}+\alpha_{k_h}=1, \ h =1,...,l,\\
&\sum_{h=1}^l\theta_h=1,\\
&\theta_h,\alpha_{1_h},\alpha_{s_h},\alpha_{k_h}\geq0, \ h =1,...,l.
\end{split}
\end{align}
\subsection{Asymptotically optimal solution}
In this section, we seek to derive the optimality conditions for (\ref{eq:5}).
We want to determine (i) the number of simulation replication allocated to each partition, (ii) the locations of the designs should be simulated in each partition and (iii) the number of simulation replications allocated to those selected designs.

According to \cite{xiao2015optimal}, we have $\theta_h/\theta_b\rightarrow0$ as the number of partitions $l$ goes to infinity.
That means the fraction of the simulation budget allocated to the partition $b$ containing the $m$-th smallest mean value far exceeds the fraction given to any other partition when $l$ goes to infinity.
Given that, the rate function $R_{m_b,i_h}(\theta_b,\theta_h, \boldsymbol{\alpha_b},\boldsymbol{\alpha_h})$ converges to $\widetilde{R}_{m_b,i_h}(\theta_h, \boldsymbol{\alpha_h})$ where
\begin{align}
\begin{split}
\widetilde{R}_{m_b,i_h}(\theta_h, \boldsymbol{\alpha_h})=\lim_{l\rightarrow \infty}R_{m_b,i_h}(\theta_b,\theta_h, \boldsymbol{\alpha_b},\boldsymbol{\alpha_h})\\
=\frac{(f(x_{m_b})-f(x_{i_h}))^2/2}{\frac{\sigma_{h}^2}{\theta_h}\Big(\frac{\eta_{i_h,1}^2}{\alpha_{1_h}}+\frac{\eta_{i_h,s}^2}{\alpha_{s_h}}+\frac{\eta_{i_h,k}^2}{\alpha_{k_h}}\Big)}.
\notag
\end{split}
\end{align}
The problem (\ref{eq:5}) can be asymptotically rewritten as
\begin{align}\label{eq:66}
\begin{split}
\max\min\Bigg\{&\!\! \min_{1\leq h\leq l\atop h\neq b}\! \bigg\{\!\min\limits_{i_h\in S_o}\!\!\widetilde{R}_{m_b,i_h}(\theta_h, \boldsymbol{\alpha_h}),\min\limits_{j_h\in S_n}\!\!\widetilde{R}_{m_b,j_h}(\theta_h, \boldsymbol{\alpha_h})\bigg\},\\
&\min \bigg\{\!\!\min\limits_{i_b\in S_o\atop  i_b\neq m_b}\!\!R_{m_b,i_b}(\theta_b,\boldsymbol{\alpha_b}),
 \min\limits_{j_b\in S_n}\!\!R_{m_b,j_b}(\theta_b,\boldsymbol{\alpha_b})\bigg\}\Bigg\}\\
&\text{s.t.}\ \alpha_{1_h}+\alpha_{s_h}+\alpha_{k_h}=1, h =1,...,l\\
&\sum_{h=1}^l\theta_h=1\\
&\theta_h,\alpha_{1_h},\alpha_{s_h},\alpha_{k_h}\geq0, h =1,...,l.
\end{split}
\end{align}
In order to better analyze problem (\ref{eq:66}) above, we decompose it as follows:
\begin{align}\label{24a}
\begin{split}
&\max\ \min_{1\leq h\leq l\atop h\neq b} \bigg\{\min\limits_{i_h\in S_o}\widetilde{R}_{m_b,i_h}(\theta_h, \boldsymbol{\alpha_h}),\min\limits_{j_h\in S_n}\widetilde{R}_{m_b,j_h}(\theta_h, \boldsymbol{\alpha_h})\bigg\}\\
&\text{s.t.}\ \alpha_{1_h}+\alpha_{s_h}+\alpha_{k_h}=1, \alpha_{1_h},\alpha_{s_h},\alpha_{k_h}\geq0,
\end{split}
\end{align}
for $h=1,2,...,l,h\neq b$, and
\begin{align}\label{23a}
\begin{split}
&\max\ \min \bigg\{\min\limits_{i_b\in S_o\atop i_b\neq m_b}R_{m_b,i_b}(\theta_b,\boldsymbol{\alpha_b}),
 \min\limits_{j_b\in S_n}R_{m_b,j_b}(\theta_b,\boldsymbol{\alpha_b})\bigg\}\\
&\text{s.t.}\ \alpha_{1_b}+\alpha_{s_b}+\alpha_{k_b}=1, \alpha_{1_b},\alpha_{s_b},\alpha_{k_b}\geq0,
\end{split}
\end{align}
for $h=b$.

\emph{Lemma 5:} \emph{Let $\theta_h^*$ and $\boldsymbol{\alpha_h^*}$ be the optimal solution to (\ref{eq:66}). The $\theta_h^*$ and $\boldsymbol{\alpha_h^*}$ are independent and can be solved separately. In addition, the $l$ $\boldsymbol{\alpha_h^*}$'s $(h=1,2,...,l)$ corresponding to the $l$ partitions are also mutually independent and can be solved separately using (\ref{24a}) and (\ref{23a}).}


Similar to Section 2, we define a key design of each partition $h$, denoted as $i_h^*$.
A false selection is most likely to happen at these key designs. That is
\begin{align}
\begin{split}
&i_h^*=\begin{cases}
\arg \!\min\limits_{1_h\leq i_h\leq k_h}\!\! \bigg\{\!\!\!\!\!\!\!&\min\limits_{i_h\in S_o}R_{m_b,i_h}(\theta_b,\theta_h, \boldsymbol{\alpha_b},\boldsymbol{\alpha_h}),\\
&\min\limits_{i_h\in S_n}R_{m_b,i_h}(\theta_b,\theta_h, \boldsymbol{\alpha_b},\boldsymbol{\alpha_h})\bigg\}, h\!\neq\! b,\\
\arg \!\min\limits_{1_b\leq i_b\leq k_{b}}\bigg\{\!\!\!\!\!\!\!\!&\min\limits_{i_b\in S_o\atop i_b\neq m_b}R_{m_b,i_b}(\theta_b,\boldsymbol{\alpha_b}),\\
&\min\limits_{i_b\in S_n}R_{m_b,i_b}(\theta_b,\boldsymbol{\alpha_b})\bigg\},h\!=\!b.\notag
\end{cases}
\end{split}
\end{align}

\subsubsection{Determine $\boldsymbol{\alpha_h^*}$ and locations of support designs}
Given Lemma 5, we can determine $\boldsymbol{\alpha_h^*}$ separately for each partition
by solving the optimization problems (\ref{24a}) and (\ref{23a}).

\begin{theorem}
The optimization problem (\ref{eq:9a}) can be asymptotically optimized with the following allocation rule:
\begin{equation}\label{eq:9}
\alpha_{r_h}^*=
\begin{cases}
\frac{|z_{r_h}|}{|z_{1_h}|+|z_{s_h}|+|z_{k_h}|}, &r_h=1_h,s_h,k_h,\\
0, &\text{otherwise},
\end{cases}
\end{equation}
where
$$
z_{r_h}^2\!=\!
\begin{cases}
\rho_{{mi}_b^*,r}^2, &h\!=\!b,  \\
\eta_{i_h^*,r}^2, &\text{otherwise},
\end{cases}
$$
for $r_h =1_h,s_h,k_h$.
\end{theorem}


Given the optimal allocation rule (Theorem 3), we next determine the optimal location for the support design $s_h$ within each partition.

\begin{theorem}
The rate function of $PFS$ with allocation rule satisfying (\ref{eq:9}) can be asymptotically maximized if the support design $s_h$ satisfies the following equations.

(a) For all $h\neq b$, support design
\begin{align}\label{eq:a12}
x_{s_h}=
\begin{cases}
x_{i_h^*}, &i_h^*\neq 1_h,k_h,\\
{(x_{1_h}+x_{k_{h}})/2}, &i_h^*=1_h,k_h,
\end{cases}
\end{align}
and the optimal allocation rule (\ref{eq:9}) becomes
\begin{align}\alpha_{r_h}^*=
\begin{cases}\label{eq:17}
1,
&r_h=i_h^*,\\
0,
&\text{otherwise}.
\end{cases}
\end{align}

(b) For $h=b$, support design
\begin{align}\label{eq:a14}
x_{s_b}\!=\!
\begin{cases}
\!x_{i^*_b}\!+\!x_{m_b}\!-\!x_{1_b},
&\frac{3x_{1_b}+x_{k_b}}{4}\leq\frac{x_{i^*_b}+x_{m_b}}{2}<\frac{x_{1_b}+x_{k_b}}{2},\\
\!x_{i^*_b}\!+\!x_{m_b}\!-\!x_{k_{b}},
&\frac{x_{1_b}+x_{k_b}}{2}<\frac{x_{i^*_b}+x_{m_b}}{2}\leq\frac{x_{1_b}+3x_{k_b}}{4},\\
\!{(x_{1_b}\!+\!x_{k_{b}})/2},
&\text{otherwise},
\end{cases}
\end{align}
and the approximate optimal allocation rule (\ref{eq:9}) is
\begin{align}
\alpha_{r_b}^*=
\begin{cases}
\frac{|\rho_{{mi}^*_b,r}|}{|\rho_{{mi}^*_b,1}|+|\rho_{{mi}^*_b,s}|+|\rho_{{mi}^*_b,k}|}, &r_b=1_b,s_b,k_b,\\
0, &\text{otherwise}.
\end{cases}
\end{align}
\end{theorem}

When the $x_{s_h}$ derived from (\ref{eq:a12}) and (\ref{eq:a14}) does not correspond to any design available, we round it to the nearest one.
\subsubsection{Determine $\theta_h^*$}
In this subsection, we aim to determine the optimal budget allocation rules between the $l$ partitions.
Since $\boldsymbol{\alpha_h^*}$ and $\theta_h^*$, can be solved separately for each partition,
$R_{m_b,i_h}(\theta_b,\theta_h, \boldsymbol{\alpha_b^*},\boldsymbol{\alpha_h^*})$ is a function of $\theta_b$ and $\theta_h$ only.
Similarly, $R_{m_b,i_h}(\theta_b, \boldsymbol{\alpha_b^*})$ is a function of $\theta_b$ only. Optimization problem (\ref{eq:5}) can be rewritten as
\begin{align}\label{eq:6}
\begin{split}
\max\min\Bigg\{\!\!\min_{1\leq h\leq l\atop h\neq b}\! \bigg\{\!\!&\min\limits_{i_h\in S_o}\!R_{m_b,i_h}(\theta_b,\theta_h, \boldsymbol{\alpha_b^*},\boldsymbol{\alpha_h^*}),\!\\
&\min\limits_{j_h\in S_n}\!R_{m_b,j_h}(\theta_b,\theta_h, \boldsymbol{\alpha_b^*},\boldsymbol{\alpha_h^*})\bigg\},\\
\min \bigg\{\!&\min\limits_{i_b\in S_o\atop i_b\neq m_b}\!R_{m_b,i_b}(\theta_b,\boldsymbol{\alpha_b^*}),
 \min\limits_{j_b\in S_n}\!R_{m_b,j_b}(\theta_b,\boldsymbol{\alpha_b^*})\bigg\}\!\!\Bigg\}\\
\text{s.t.}\ \sum_{h=1}^l&\theta_h=1,\theta_h\geq0, h =1,...,l.
\end{split}
\end{align}
By solving the optimization model (\ref{eq:6}), we can get the optimal $\theta_h^*$ for all $h =1,...,l$.
\begin{theorem}
The optimal allocation that asymptotically minimizes the $PFS$ for the problem (\ref{eq:9a}) is that
$
\theta_h^*=\frac{\gamma_h^*}{\sum\limits_{i=1}^l\gamma_i^*},
$
where
${\gamma_b^*}=\sigma_{b}\sqrt{\Big(\frac{\eta_{m_b,1}^2}{\alpha^*_{1_b}}\!+\!\frac{\eta_{m_b,s}^2}{\alpha^*_{s_b}}\!+\!\frac{\eta_{m_b,k}^2}{\alpha^*_{k_b}}\!\Big)\sum\limits_{h\neq b}\frac{{\gamma_h^*}^2}{\sigma_{h}^2}}$ and $\gamma_h^*=\sigma_h^2/(f(x_{m_b})-f(x_{i_h^*}))^2,\ h\neq b.
$

\end{theorem}

The results of Theorem 5 indicate the optimal simulation budget allocation between the $l$ partitions. The results of Theorem 4 determine which designs should be selected for simulation and the number of simulation replications allocated to these selected designs. Since the solution space is divided into adjacent partitions, this optimal budget allocation procedure can be used for more general cases where the underlying function is non-quadratic. By conducting simulation replications on fewer designs, the simulation efficiency is dramatically enhanced compared to the existing methods.

\section{Sequential Budget Allocation Algorithm}
Based on the discussion above, we propose two sequential budget allocation algorithm to implement the optimality conditions. They are Optimal Optimal Computing Budget Allocation for selecting the top-$m$ designs with regression (OCBA-mr) based on Theorems 1 and 2 and Optimal Optimal Computing Budget Allocation for selecting the top-$m$ designs with regression in partitioned domains (OCBA-mrp) based on Theorems 4 and 5.

\noindent\textbf{(1) OCBA-mr Algorithm}

\noindent\textbf{INITIALIZE:} $\kappa\leftarrow 0$; Perform $n_0$ simulation replications for three design locations in each partition; by convention we use the D-optimal support designs, i.e., $n_{1}^\kappa=n_{{((1+t)/2)}}^\kappa=n_{t}^\kappa=n_0$ (round as needed). The incremental simulation budget is $\Delta$.

\noindent\textbf{LOOP: WHILE} $\sum_{i=1}^tn_{i}^\kappa<n$ $\textbf{DO}$

\noindent\textbf{UPDATE:} Estimate a quadratic regression equation $\hat{f}(x_{i})$ using the least squares method based on the information from all prior simulation runs.

Calculate the mean and variance of each design using $\hat{f}(x_{i})=\hat{\beta}_{0}+\hat{\beta}_{1}x_{i}+\hat{\beta}_{2}x_{i}^2.$

Determine the observed subset $\hat{S}_o$ and $\hat{S}_n$; the design with the $m$-th smallest sample mean value; the observed key design $\hat{i}^*$.

Calculate $\boldsymbol{\hat\alpha^*}$ and determine the support designs using Theorems 1 and 2.

\noindent\textbf{ALLOCATE:} Increase the simulation budget by $\Delta$ and calculate the new budget allocation rule using $\boldsymbol{\hat\alpha^*}$.

\noindent\textbf{SIMULATE:} Perform $\max(0, n_{i}^\kappa-n_{i}^{\kappa-1})$ simulation replications for design ${i}$ and $i=1,s,t$. $\kappa\leftarrow \kappa+1$.

\noindent\textbf{END OF LOOP.}

\noindent\textbf{(2) OCBA-mrp Algorithm}

\noindent\textbf{INITIALIZE:} $\kappa\leftarrow 0$; Perform $n_0$ simulation replications for three design locations in each partition; by convention we use the D-optimal support designs, i.e., $n_{1_h}^\kappa=n_{{((1+k^{(h)})/2)}_h}^\kappa=n_{k_h}^\kappa=n_0$ (round as needed). The incremental simulation budget is $\Delta$.

\noindent\textbf{LOOP: WHILE} $\sum_{h=1}^l\sum_{i_h=1_h}^{k_{h}}n_{i_h}^\kappa<n$ $\textbf{DO}$

\noindent\textbf{UPDATE:} Estimate the quadratic regression equations $\hat{f}(x_{i_h})$ for each partition using the least squares method based on the information from all prior simulation runs.

Calculate the means and variances of each design using $\hat{f}(x_{i_h})=\hat{\beta}_{h0}+\hat{\beta}_{h1}x_{i_h}+\hat{\beta}_{h2}x_{i_h}^2,\ i_h=1_h,2_h,...,k_h, h=1,2,...,l.$

Determine the observed subset $\hat{S}_o$ and $\hat{S}_n$; the design with the $m$-th smallest sample mean value; the observed key design $\hat{i_h}^*$ of each partition.

Calculate $\boldsymbol{\hat\alpha_h^*}$ and determine the support designs in each partition using Theorem 4.
Calculate ${\hat\theta_h^*}$ using Theorem 5.

\noindent\textbf{ALLOCATE:} Increase the simulation budget by $\Delta$ and calculate the new budget allocation rule using $\boldsymbol{\hat\alpha_h^*}$ and $\hat\theta_h^*$.

\noindent\textbf{SIMULATE:} Perform $\max(0, n_{i_h}^\kappa-n_{i_h}^{\kappa-1})$ simulation replications for design ${i_h}$ in each partition $h=1,2,...,l$ and $i_h=1_h,s_h,k_h$. $\kappa\leftarrow \kappa+1$.

\noindent\textbf{END OF LOOP.}

\section{Numerical Experiments}
In this section, we test our proposed simulation budget allocation rules, OCBA-mr and OCBA-mrp, with some existing methods on several typical optimal subset selection problems.
We use the following three budget allocation approaches for comparison. The top-$m$ designs are selected based on their mean values.
\begin{itemize}
    \item \emph{Equal Allocation (EA)}: The EA is the most commonly used and simplest method, which allocates the simulation budget equally to each of the design.

	\item \emph{OCBAm+ Allocation}: The OCBAm+ is an efficient simulation budget allocation procedure for selecting the top $m$ designs. It was developed based on the OCBA framework. For detail see \cite{zhang2016simulation}.

	\item \emph{OCBA-ss Allocation}: The OCAB-ss procedure was proposed in \cite{gao2016new} and seek to solve the optimal subset R\&S problems for general underlying distributions using the large deviations theory.

     \item \emph{OSD Allocation}: The OSD procedure was proposed in \cite{brantley2013efficient} and seek to solve the single best R\&S problems using regression metamodel.

\end{itemize}

The OCBA-mrp, OCBA-mr and OSD estimate the performance of each design based on some regression metamodels.
The EA, OCBAm+ and OCBA-ss use the mean value of each design for comparison, and the mean performance value is computed directly from the simulation output. Unlike our proposed methods and OSD, they do not rely on any response surface to estimate the performance value for any design.
In order to compare the performance of these allocation approaches, we test them empirically on the following experiments.

\begin{itemize}

  \item Experiment 1: Consider the optimization problem
  \begin{align}
  \begin{split}
  &\min\ f(x)=(x-5)^2.\\
  \end{split}\notag
  \end{align}
  The experiment assumes that the noise for simulation has normal distribution $N(0,2^2)$
  for the objective value. We discretize the domain of the function into 100 evenly spaced points from 0 to 10. Since the underlying function is quadratic, we can use the OCBA-mr method directly for this problem. In order to utilize the OCBA-mrp method, we divide the 100 points into 5 adjacent partitions and there are 20 points in each partition. We want to select the top $m=5$ designs.

  \item Experiment 2: Consider the Griewank function
  \begin{align}
  \begin{split}
  \min\ f(x)=10\times(1+(1/4000)\times x^2-\cos(x)).\notag
  \end{split}
  \end{align}
  We discretize the domain of the function into 100 evenly spaced points from 0 to 20. The nature of this underlying function does not allow us to apply the OCBA-mr method directly. To test the performance of the OCBA-mr method for the non-quadratic problem, we make a slight modification by dividing the solution space into 5 adjacent partitions. Then, we allocate the simulation budget equally to each partition and use the OCBA-mr method within each partition.
  In this experiment, the OCBA-mrp is also tested under this partition pattern.
  The experiment assumes that the noise for simulation has normal distribution $N(0, 0.2^2)$
  We want to select the top $m=3$ designs.


  \item Experiment 3: Consider the optimization problem
  \begin{align}
  \begin{split}
  \min f(x)\!=\!\sin(x)\!+\!\sin(10x/3)\!+\!\log(x)\!-\!0.84x\!+\!3.\notag
  \end{split}
  \end{align}
 The entire domain consists of 200 discrete points for $x\in[0, 8]$. The 200 points are further divided into 10 adjacent partitions so that there are 20 points in each partition.
  This experiment assumes that the simulation noise is a standard normal random variable.
  Similar to the setting in experiment 2, the OCBA-rm method is modified to handle this non-quadratic problem.
  We want to select the top $m=5$ designs.

  \item Experiment 4: Consider the optimization problem
  \begin{align}
  \begin{split}
  \min\ f(x)=2(x-0.75)^2+\sin(8\pi x-\pi/2).\notag
  \end{split}
  \end{align}
  The experiment assumes that the noise for simulation has standard normal distribution
  for the objective measure. We discretize the domain of the function into 200 evenly spaced points from 0 to 2. The 200 points are further divided into 20 adjacent partitions so that there are 10 points in each partition. Similar to the setting in experiment 2, the OCBA-rm method is modified to handle this non-quadratic problem.
  We want to select the top $m=3$ designs.

  \item Experiment 5: Consider the optimization problem
  \begin{align}
  \begin{split}
  \min\ f(x_1,x_2)=\frac{1}{40}(x_1^2+x_2^2)-\cos(x_1)\cos(\frac{x_2}{\sqrt{2}})+1.\notag
  \end{split}
  \end{align}
  The experiment assumes that the noise for simulation has normal distribution $N(0,2^2)$
  for the objective measure. $x_1$ and $x_2$ are continuous variables with $-5\leq x_1\leq 5$ and $-5\leq x_2\leq 5$. We discretize the solution space into $11\times 11$ discrete points $x_1=\{-5,-4,...,5\}$ and $x_2=\{-5,-4,...,5\}$.
  We divide the $11\times 11$ designs into 11 adjacent partitions with 11 designs in each partition. For each $i=1,2,...,11$, partition $i$ consists of design points $(-5, i-6)$, $(-4, i-6)$, ..., $(5, i-6)$.
  We want to select the top $m=3$ designs, which are $(0,-1)$, $(0,0)$ and $(0,1)$.

  \item Experiment 6: Consider the $(s,S)$ inventory problem in the Simulation Optimization Library
      (\url{http://simopt.org/wiki/index.php?title=SS_Inventory}). The decision variable
      $(s,S)$ is the inventory strategy. When the inventory of design $(s,S)$ on hand below $s$, an order is incurred to supplement the inventory to $S$. The optimization problem is to minimize the the E[Total cost per period]. The domain is discretized into $20\times 20$ discrete points with $s=\{810,820,...,1000\}$ and $S=\{1510,1520,...,1700\}$. The discrete points are divided into 20 partitions with 20 points in each partition. For each $i=1, 2,...,20$, partition $i$ consists of points $(800+10i,1510),(800+10i,1520),...,(800+10i,1700)$. In each partition, we assume the performances of the E[Total cost per period] of each point can be modeled as quadratic functions, where the independent variable is $S$ and $s$ is kept as a constant. We want to select the top $m=3$ designs, which are $(810,1510)$, $(810,1520)$ and $(820,1510)$.


\end{itemize}

For the budget allocation procedures mentioned above, we let the initial number of replication $n_0$ be 10 and the incremental budget $\Delta$ be 100 for all the five experiments. Figures 1 reports the comparison results in terms of $PCS$. The estimate of $PCS$ is based on the average of 2000 independent replications of each procedure for experiments 1-6.


\begin{figure*}[htbp]
\begin{center}
\subfigure[Experiment 1]{
\includegraphics[height=1.7in,width=2.2in]{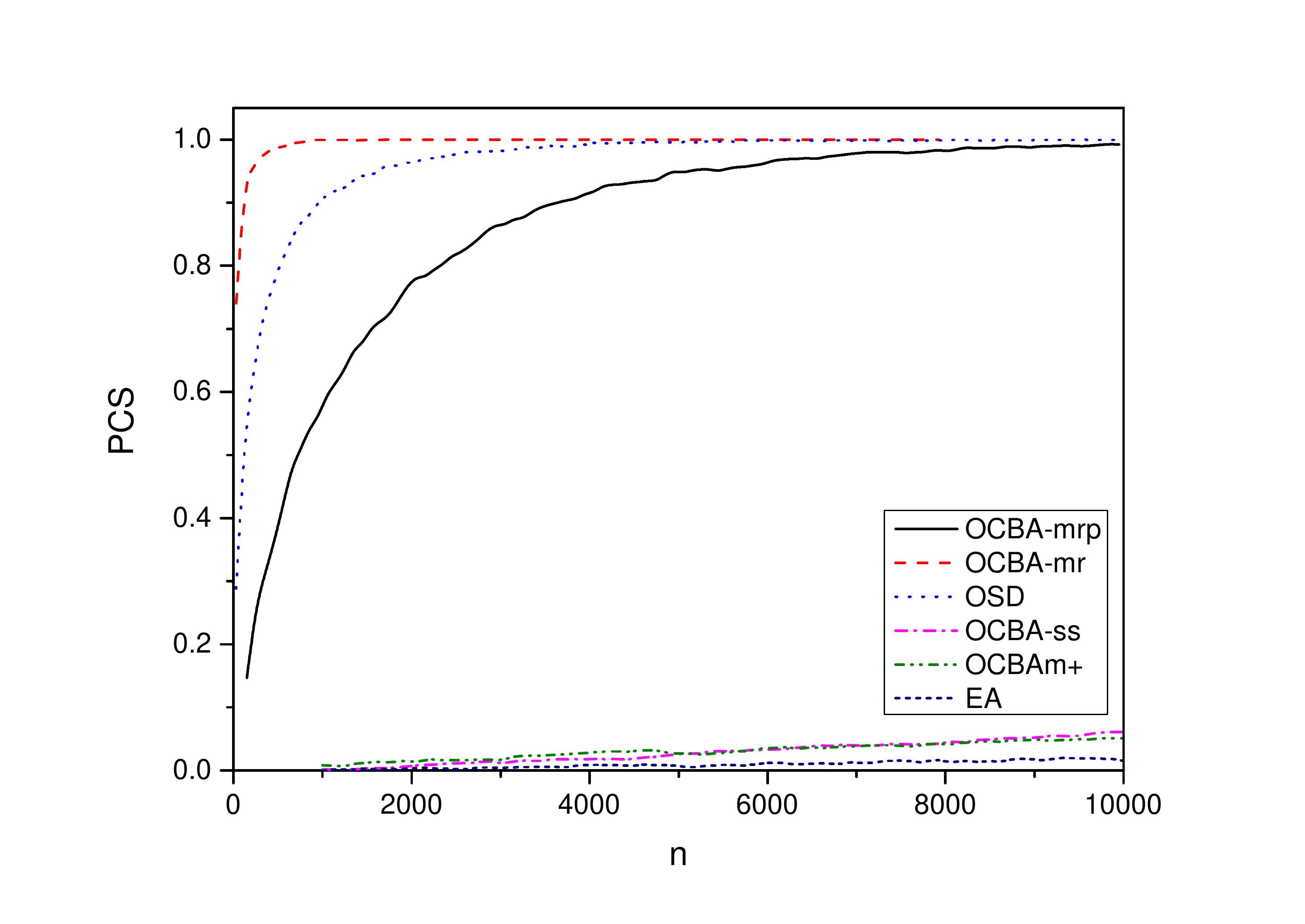}
\label{fig:ex1}}
\subfigure[Experiment 2]{
\includegraphics[height=1.7in,width=2.2in]{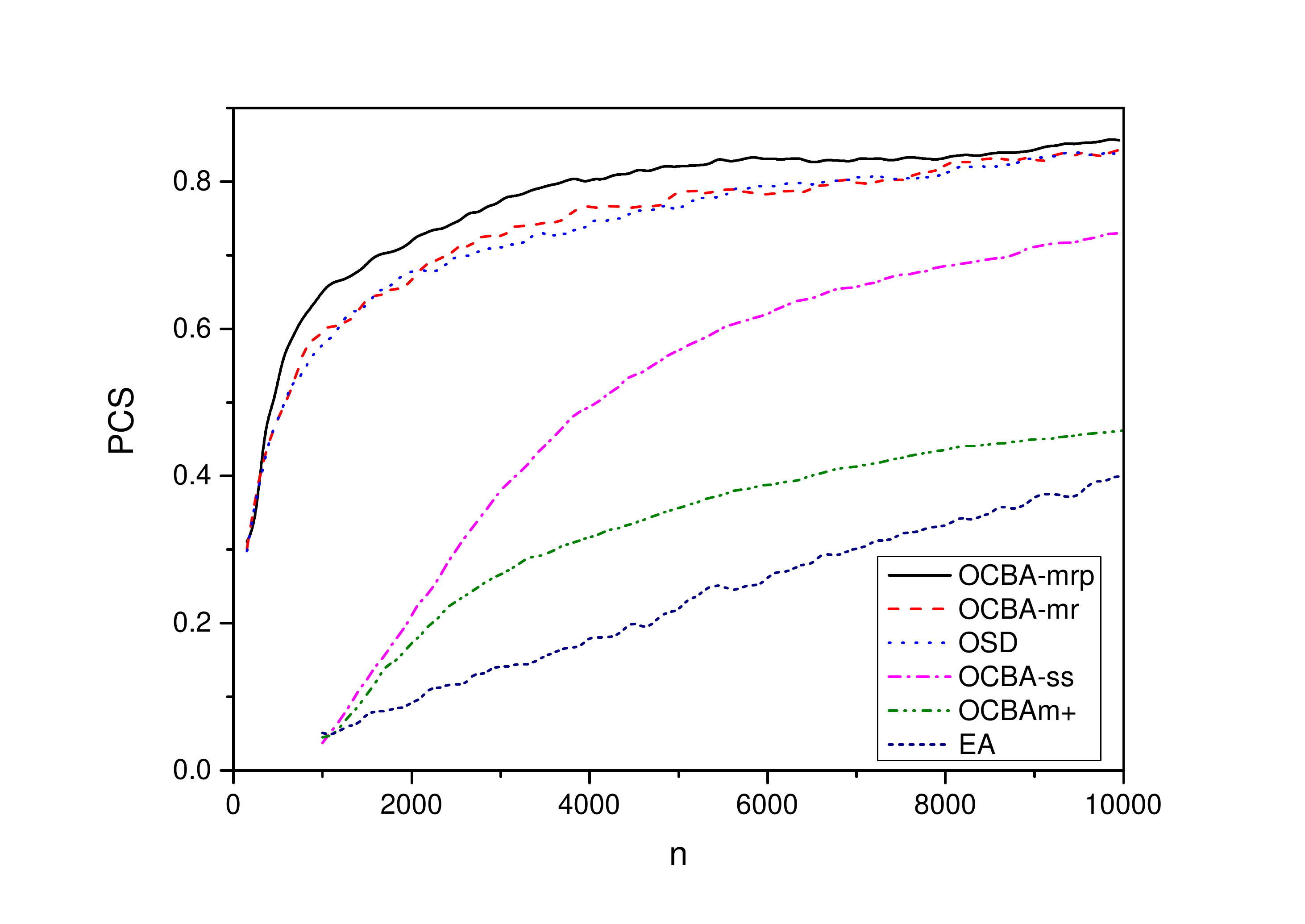}
\label{fig:ex2}}
\subfigure[Experiment 3]{
\includegraphics[height=1.7in,width=2.2in]{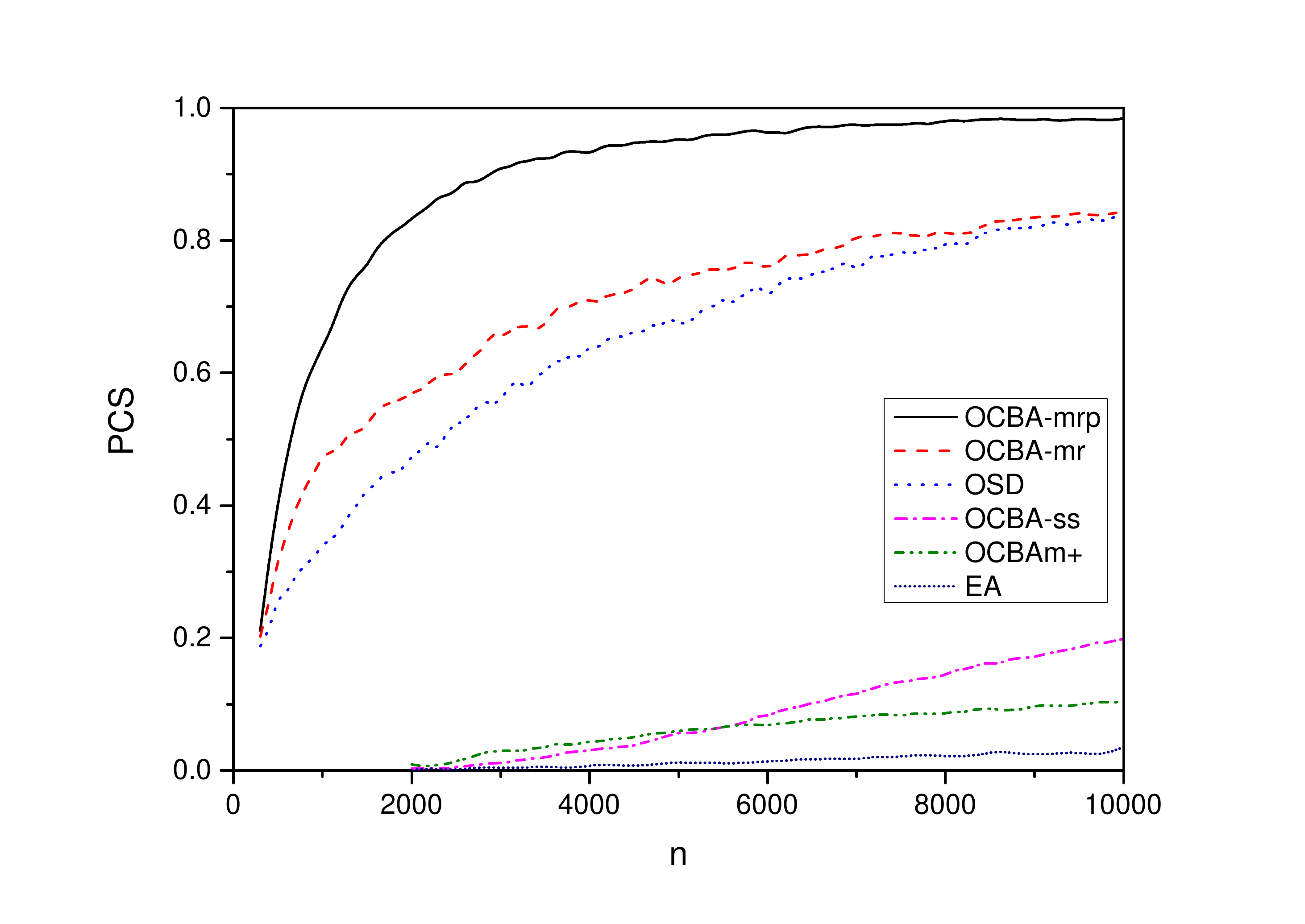}
\label{fig:ex3}}
\subfigure[Experiment 4]{
\includegraphics[height=1.7in,width=2.2in]{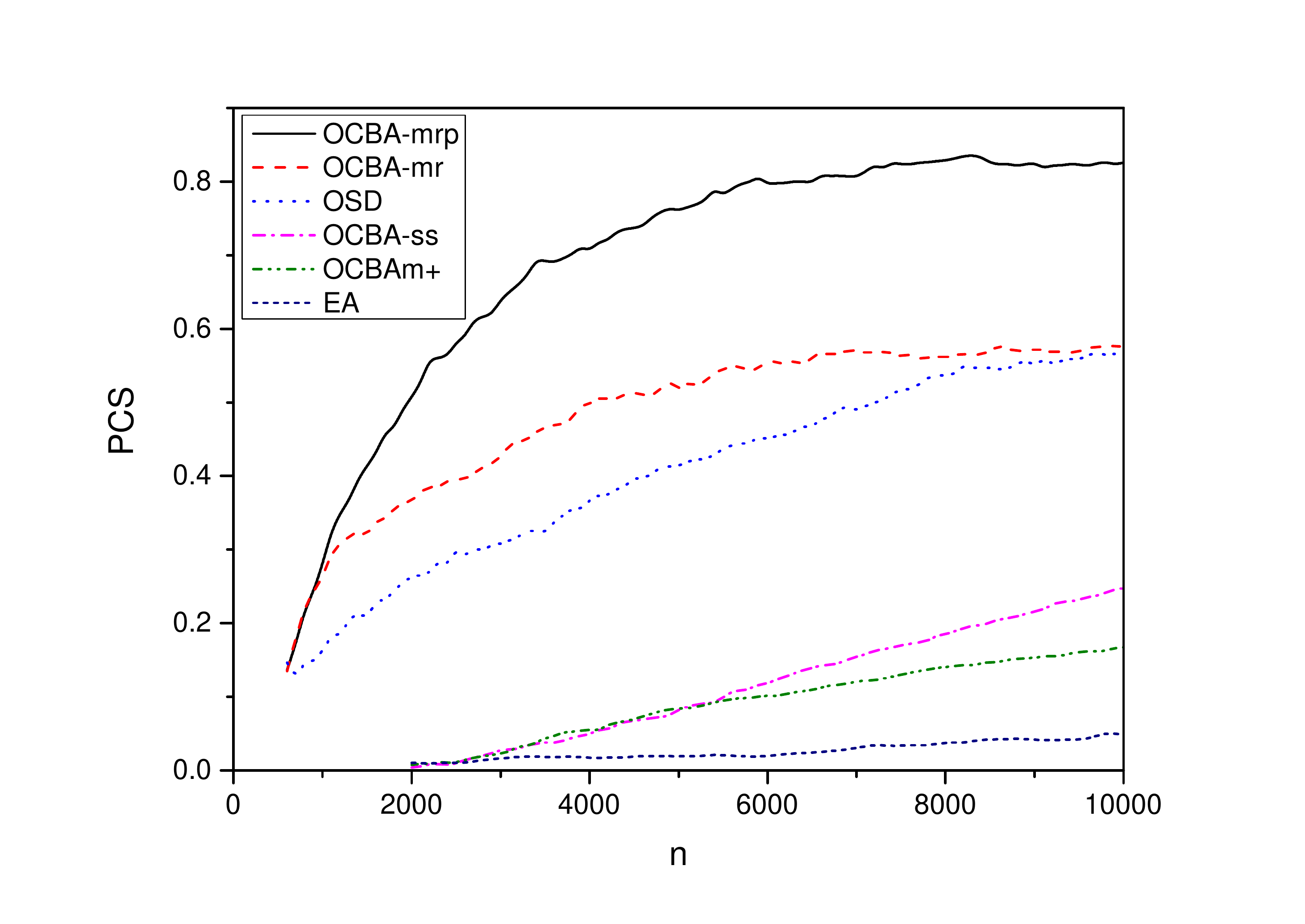}
\label{fig:ex4}}
\subfigure[Experiment 5]{
\includegraphics[height=1.7in,width=2.2in]{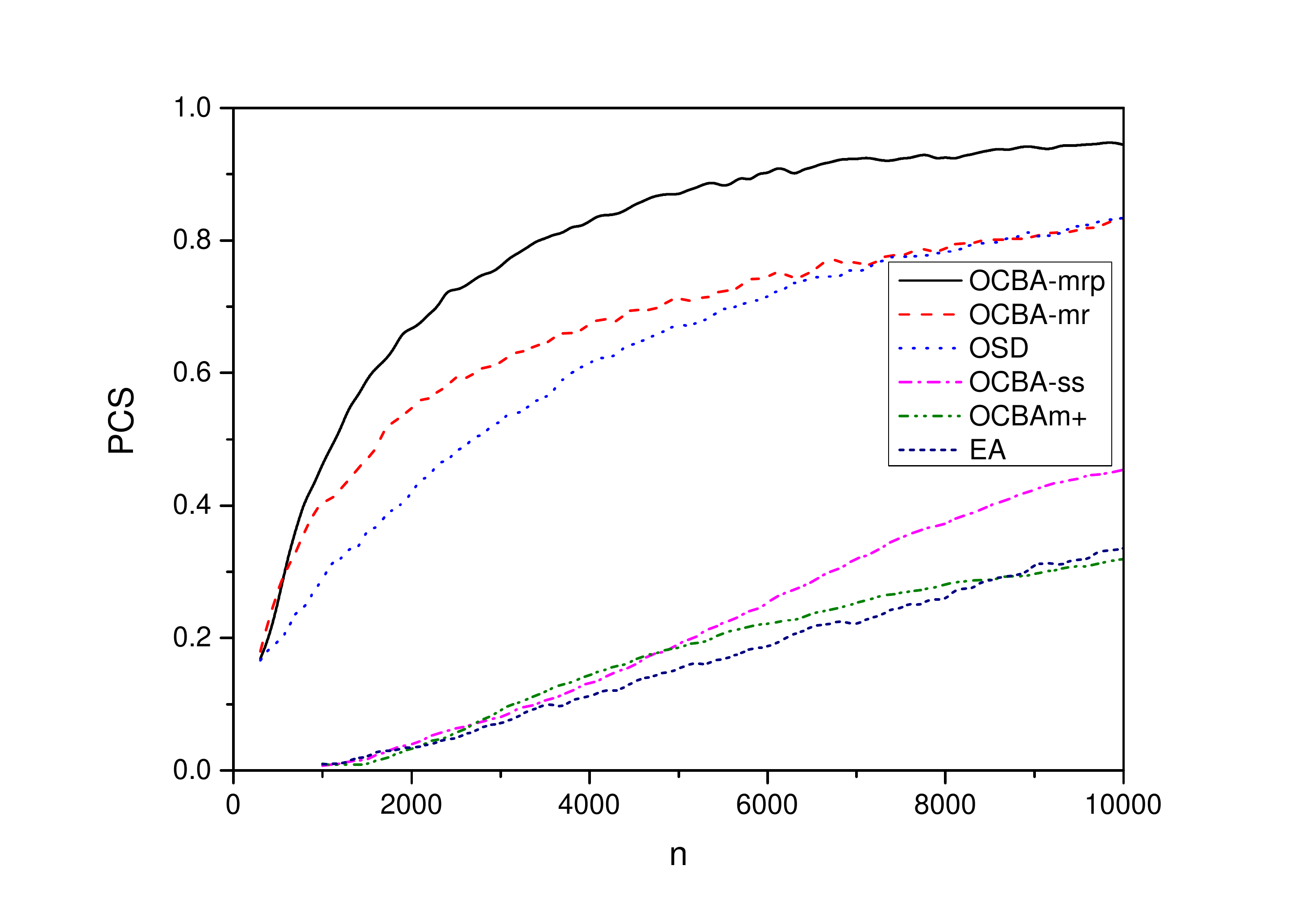}
\label{fig:ex5}}
\subfigure[Experiment 6]{
\includegraphics[height=1.7in,width=2.2in]{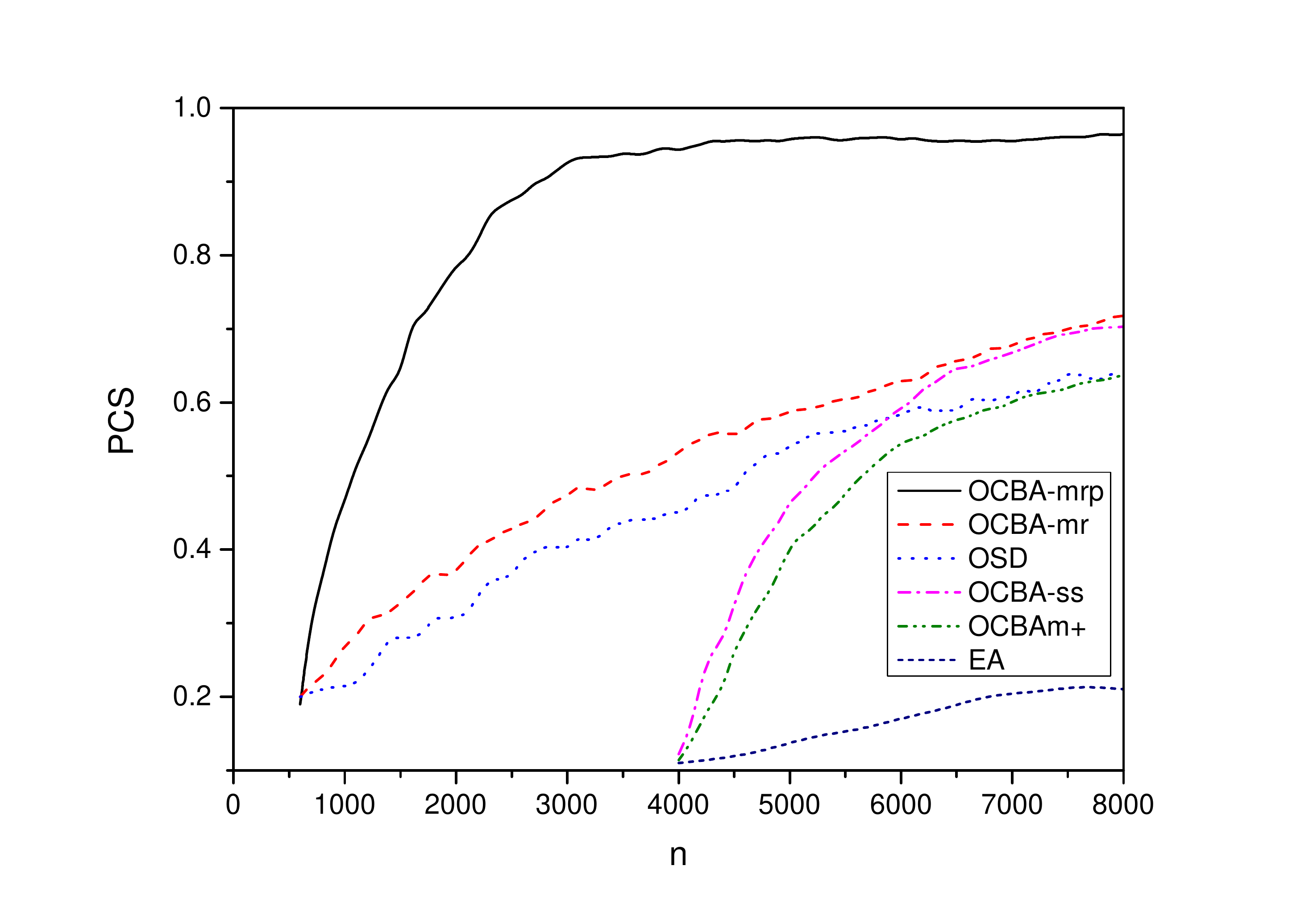}
\label{fig:ex6}}
\caption{Comparison results for different test examples.}
\label{fig:algorithm}
\end{center}
\end{figure*}

It is observed that all the procedures obtain higher $PCS$ as the available simulation budget increases. Our proposed procedures OCBA-mr and OCBA-mrp perform the best on the tested examples.
It shows that our proposed procedures have dramatically improved the selection quality of the top-$m$ designs compared to the existing methods.
When the considered problem is quadratic or approximately quadratic, such as example 1, the OCBA-mr method has the best performance. That is because the OCBA-mr method is developed for the quadratic or approximately quadratic problem setting. It can make full use of its advantages by simulating fewer designs once the quadratic assumption is satisfied.
When the considered problem is partially quadratic or non-quadratic, such as examples 2-6, the OCBA-mrp method performs better than the modified OCBA-mr method. Because the OCBA-mrp can intelligently allocate the simulation budget both between and within each partition.


Compared to the existing approaches, our proposed approaches incorporate the information from the simulated designs into regression metamodel(s). We only need to simulate a subset of the alternative designs to build up the regression metamodel(s). The performances of the designs that have not been simulated can be inferred based on the metamodel(s). It dramatically improves the selection efficiency.

\section{Conclusions}
In this study, we further enhanced the simulation efficiency of finding the top-$m$ designs by incorporating quadratic regression equation(s). Using the large deviations theory, we formulate the problem of selecting the top-$m$ designs as that of maximizing the minimum convergence rate of the false selection probability. We derived two asymptotically optimal budget allocation rules, OCBA-mr and OCBA-mrp, for one partition and multi-partition problem setting, respectively.
Numerical results suggest that the proposed approaches can dramatically improve the selection efficiency compared to the existing R\&S approaches. When the underlying function can be approximated by a quadratic function across the solution space, the OCBA-mr is more efficient. For the non-quadratic problem setting, the OCBA-mrp demonstrates the best performance.

\bibliographystyle{apalike}
\bibliography{library}

\newpage

\appendix{Online Supplement}

\section*{Proof of Lemma 2}

Let $\Lambda_{i}^{(d)}(\theta)=\ln E(e^{\theta \hat{d}_{i,m}})$ denote the log-moment generating function of $\hat{d}_{i,m}$ and $I_{i}^{(d)}(\cdot)$ denote the Fenchel-Legendre transform of $\Lambda_{i}^{(d)}$, that is
$$I_{i}^{(d)}(x)=\sup_{\theta\in \mathbb{R}}(\theta x-\Lambda_{i}^{(d)}(\theta)).$$

For any set $A$, let $A^\circ$ denotes its interior and for any function $f(\cdot)$, let $f'(x)$ denote the derivative of $f$ at $x$. The effective domain of $\Lambda_{i}$ is $D_{\Lambda_{i}}=\{\theta\in\mathbb{R}:\Lambda_{i}(\theta)<\infty\}$ and $F_{i}=\{\Lambda_{i}'(\theta):\theta\in D_{\Lambda_{i}}^\circ\}$. Let $f_{min}=\min_{i\in \{1,2,...,t\}}f(x_{i})$ and $f_{max}=\max_{i\in \{1,2,...,t\}}f(x_{i})$. We assume that the interval $[f_{min},f_{max}]\subset\bigcap_{i=1}^t F_{i}^\circ$. This assumption can be satisfied by some common families of distributions such Normal, Bernoulli, Poisson and Gamma family \citep{glynn2004large}. It ensures that the range of $\hat{f}(x_{i})$
include $[f_{min},f_{max}]$ and $P(\hat{f}(x_{i})>\hat{f}(x_{m}))$ and $P(\hat{f}(x_{j})<\hat{f}(x_{m}))$ are positive for $i\in S_o$ and $j\in S_n$.

Then, $\hat{d}_{i,m}=\hat{f}(x_{i})-\hat{f}(x_{m})$ satisfies the large deviation principle with good rate functions
\begin{equation}
\begin{split}
&\left.
\begin{aligned}
\!-\!\lim_{n\rightarrow \infty}\frac{1}{n}\log P\{\hat{d}_{i,m}\!>\! 0\},\ &i\!\in\! S_o, i\!\neq\! m\\
\!-\!\lim_{n\rightarrow \infty}\frac{1}{n}\log P\{\hat{d}_{i,m}\!<\! 0\},\ &i\!\in\! S_n
\end{aligned}\right\}\!=\!R_{m,i}(\boldsymbol{\alpha})
 \!=\!I_{i}^{(d)}(0),
\end{split}\notag
\end{equation}

Since $\hat{d}_{i,m}=\hat{f}(x_{i})-\hat{f}(x_{m})$ are all normally distributed random variables, the rate function can be expressed as follow according to the results in \cite{glynn2004large}.
$$R_{m,i}(\boldsymbol{\alpha})
 = \frac{({f}(x_{m}) - {f}(x_{i}))^2/2}{{\sigma^2} \Big( \frac{\rho_{mi ,1}^2}{\alpha_{1}} + \frac{\rho_{mi ,s}^2}{\alpha_{s}} + \frac{\rho_{mi ,t}^2}{\alpha_{t}} \Big)}.$$

\section*{Proof of Lemma 3}

According to the definition of $PCS$, we have
$$\begin{aligned}
&PCS=P\bigg\{\bigcap_{i\in S_o}\bigcap_{j\in S_n}\hat{f}(x_{i})\leq\hat{f}(x_{j})\bigg\}\\
&\geq P\bigg\{\Big(\bigcap_{i\in S_o\atop i\neq m}\hat{f}(x_{i})\leq \hat{f}(x_{m})\Big)\bigcap
\Big(\bigcap_{j\in S_n}\hat{f}(x_{j})\geq \hat{f}(x_{m})\Big)\bigg\}\\
&\geq1-\sum_{i\in S_o,\atop i\neq m}P\{\hat{f}(x_{i})\geq \hat{f}(x_{m})\}
-\sum_{j\in S_n}P\{\hat{f}(x_{j})\leq \hat{f}(x_{m})\}.
\end{aligned}$$
The last inequality follows from the Bonferroni inequality.
Since $PFS=1-PCS$, we can get that
\begin{align}
PFS\leq\sum_{i\in S_o\atop i\neq m}P\{\hat{f}(x_{i})\geq \hat{f}(x_{m})\}
+\sum_{j\in S_n}P\{\hat{f}(x_{j})\leq \hat{f}(x_{m})\}.\notag
\end{align}
Therefore, a lower bound on $PFS$ can be derived as
$$
\begin{aligned}
\max \bigg\{\!\max\limits_{i\in S_o\atop i\!\neq\! m}P\{\hat{f}(x_{i})\!\geq\! \hat{f}(x_{m})\},
 \max\limits_{j\in S_n}P\{\hat{f}(x_{j})\!\leq\! \hat{f}(x_{m})\}\bigg\}
\end{aligned}$$
and a upper bound on $PFS$ can be derived as
$$
\begin{aligned}
|S_o|\!\times\!|S_n|\!\times\!\max \bigg\{\!\max\limits_{i\in S_o\atop i\neq m}P\{\hat{f}(x_{i})\!\geq\! \hat{f}(x_{m})\},
 \max\limits_{j\in S_n}P\{\hat{f}(x_{j})\!\leq\! \hat{f}(x_{m})\}\!\bigg\}.
\end{aligned}
$$
According to Lemma 2, as $n$ goes to infinity we can get the convergence rate function shown in Lemma 3.

\section*{Proof of Theorem 1}

Given the definition of the key design $i^*$, we rewrite (\ref{eq:3a}) as
\begin{align}\label{eq:5a}
\begin{split}
&\max\ {R}_{m,i^*}(\boldsymbol{\alpha})\\
\text{s.t.}\
&\alpha_{1}+\alpha_{s}+\alpha_{t}=1,
\alpha_{1},\alpha_{s},\alpha_{t}>0.
\end{split}
\end{align}

Since ${R}_{m,i^*}(\boldsymbol{\alpha})$ is concave and strictly increasing functions of $\boldsymbol{\alpha}$ \citep{xiao2015optimal,glynn2004large},
the optimization problem is a convex programming problem. 
We define a Lagrangian function $L=-{R}_{m,i^*}( \boldsymbol{\alpha})+\lambda(\alpha_{1}+\alpha_{s}+\alpha_{t}-1)$ and determine the optimal allocation based on its Karush-Kuhn-Tucker (KKT) conditions.
\begin{align}\label{eq6}
\begin{split}
\frac{\partial L}{\partial\alpha_{r}}=-\frac{(f(x_{m})-f(x_{i^*}))^2/2}
{{\sigma^2}\Big(\frac{\rho_{mi^*,1}^2}{\alpha_{1}}+\frac{\rho_{mi^*,s}^2}{\alpha_{s}}+\frac{\rho_{mi^*,t}^2}{\alpha_{t}}\Big)^2}\frac{\rho_{mi^*,r}^2}{{\alpha_{r}^*}^2}\!+\!\lambda\!=\!0,\\
{r}=1,s,t.
\end{split}
\end{align}
Plug (\ref{eq6}) into $\alpha_{1}^*+\alpha_{s}^*+\alpha_{t}^*=1$, we can establish
\begin{equation}\label{24b}
\alpha_{r}^*=\frac{|\rho_{mi^*,r}|}{|\rho_{mi^*,1}|+|\rho_{mi^*,s}|+|\rho_{mi^*,t}|}, {r}=1,s,t.
\end{equation}
Then, we can get the conclusions in Theorem 1.

\section*{Proof of Theorem 2}

Plugging (\ref{eq:4a}) into ${R}_{m,i^*}(\boldsymbol{\alpha})$, we have
$$\!{R}_{m,i^*}(\boldsymbol{\alpha})=\frac{(f(x_{m})-f(x_{i^*}))^2/2}
{{\sigma^2}\left(|\rho_{mi^*,1}|+|\rho_{mi^*,s}|+|\rho_{mi^*,t}|\right)^2}.$$
\begin{align}
\begin{split}
\frac{\partial {R}_{m,i^*}(\boldsymbol{\alpha})}{\partial x_{s}}=&\frac{-(f(x_{m})-f(x_{i^*})^2}
{{\sigma^2}\Big(|\rho_{mi^*,1}|+|\rho_{mi^*,s}|+|\rho_{mi^*,t}|\Big)^3}\\
&\times\Big(\frac{\partial|\rho_{mi^*,1}|}{\partial x_{s}}\!+\!\frac{\partial|\rho_{mi^*,s}|}{\partial x_{s}}\!+\!\frac{\partial|\rho_{mi^*,t}|}{\partial x_{s}}\Big).
\notag
\end{split}
\end{align}
where
$$\frac{\partial \rho_{mi^*,1}}{\partial x_{s}}=\frac{(x_{i^*}-x_{m})(x_{i^*}+x_{m}-x_{1}-x_{t})}{(x_{1}-x_{s})^2(x_{1}-x_{t})},$$
\begin{align}
\begin{split}
\frac{\partial \rho_{mi^*,s}}{\partial x_{s}}
=-\frac{(x_{i^*}\!-\!x_{m})(x_{i^*}\!+\!x_{m}\!-\!x_{1}\!-\!x_{t})(2x_{s}\!-\!x_{1}\!-\!x_{t})}{(x_{s}-x_{1})^2(x_{s}-x_{t})^2},
\notag
\end{split}
\end{align}
$$\frac{\partial \rho_{mi^*,t}}{\partial x_{s}}=\frac{(x_{i^*}-x_{m})(x_{i^*}+x_{m}-x_{1}-x_{t})}{(x_{t}-x_{1})(x_{t}-x_{s})^2}.$$

In order to analyze the optimal location for $x_{s}$, we consider the following five cases.

Case I: $\frac{x_{i^*}+x_{m}}{2}<\frac{3x_{1}+x_{t}}{4}$

We have ${\partial {R}_{m,i^*}(\boldsymbol{\alpha})}/{\partial x_{s}}>0$ when $x_s<{(x_{1}+x_{t})/2}$; ${\partial {R}_{m,i^*}(\boldsymbol{\alpha})}/{\partial x_{s}}<0$ when $x_s>{(x_{1}+x_{t})/2}$. Therefore, when ${(x_{i^*}+x_{m})}/{2}<{(3x_{1}+x_{t})}/{4}$, we choose $x_s={(x_{1}+x_{t})/2}$.

Case II: $\frac{3x_{1}+x_{t}}{4}\leq\frac{x_{i^*}+x_{m}}{2}<\frac{x_{1}+x_{t}}{2}$

We have ${\partial {R}_{m,i^*}(\boldsymbol{\alpha})}/{\partial x_{s}}>0$ when $x_s<x_{i^*}+x_{m}-x_{1}$; ${\partial {R}_{m,i^*}(\boldsymbol{\alpha})}/{\partial x_{s}}<0$ when $x_s>x_{i^*}+x_{m}-x_{1}$. Therefore, when ${(3x_{1}+x_{t})}/{4}\leq{(x_{i^*}+x_{m})}/{2}<{(x_{1}+x_{t})}/{2}$, we choose $x_s=x_{i^*}+x_{m}-x_{1}$.

Case III:  $\frac{x_{i^*}+x_{m}}{2}>\frac{x_{1}+3x_{t}}{4}$

We have ${\partial {R}_{m,i^*}(\boldsymbol{\alpha})}/{\partial x_{s}}>0$ when $x_s<{(x_{1}+x_{t})/2}$; ${\partial {R}_{m,i^*}(\boldsymbol{\alpha})}/{\partial x_{s}}<0$ when $x_s>{(x_{1}+x_{t})/2}$. Therefore, when ${(x_{i^*}+x_{m})}/{2}>{(x_{1}+3x_{t})}/{4}$, we choose $x_s={(x_{1}+x_{t})/2}$.

Case IV: $\frac{x_1+x_t}{2}<\frac{x_{i^*}+x_m}{2}\leq\frac{x_1+3x_t}{4}$

We have ${\partial {R}_{m,i^*}(\boldsymbol{\alpha})}/{\partial x_{s}}>0$ when $x_s<x_{i^*}+x_{m}-x_{t}$; ${\partial {R}_{m,i^*}(\boldsymbol{\alpha})}/{\partial x_{s}}<0$ when $x_s>x_{i^*}+x_{m}-x_{t}$. Therefore, when ${(3x_{1}+x_{t})}/{4}\leq{(x_{i^*}+x_{m})}/{2}<{(x_{1}+x_{t})}/{2}$, we choose $x_s=x_{i^*}+x_{m}-\!x_{t}$.

Case V: $\frac{x_{i^*}+x_{m}}{2}=\frac{x_{1}+x_{t}}{2}$

The location of $x_s$ has no influence on the results. Therefore, we choose $x_s={(x_{1}+x_{t})/2}$ for consistency with Case II and Case III.

Analyzing the above the results, we can get the conclusions in Theorem 2.

\section*{Proof of Lemma 4}

Let $\Lambda_{i_h}(\theta)=\ln E(e^{\theta \hat{f}(x_{i_h})})$ denote the log-moment generating function of $\hat{f}(x_{i_h})$ and $I_{i_h}(\cdot)$ denote the Fenchel-Legendre transform of $\Lambda_{i_h}$, that is
$$I_{i_h}(x)=\sup_{\theta\in \mathbb{R}}(\theta x-\Lambda_{i_h}(\theta)).$$
Let the scaled cumulant generating function of $(\hat{f}(x_{m_b}),\hat{f}(x_{i_h}))$ be denoted as follows:

\begin{equation}
\begin{split}
\lim_{n\rightarrow\infty}\frac{1}{n}\Lambda^{(\hat{f}(x_{m_b}),\hat{f}(x_{i_h}))}(n\theta_b,n\theta_h)\\
=\lim_{n\rightarrow\infty}\frac{1}{n}\ln E(e^{n\theta_b\hat{f}(x_{m_b})+n\theta_h\hat{f}(x_{i_h})}),
\notag
\end{split}
\end{equation}
By the G\"{a}rtner-Ellis theorem, $(\hat{f}(x_{m_b}),\hat{f}(x_{i_h}))$ satisfy the large deviation principle with good rate functions
\begin{equation}
\begin{split}
&\left.
\begin{aligned}
-\lim_{n\rightarrow \infty}\frac{1}{n}\log P\{\hat{f}(x_{i_h})> \hat{f}(x_{m_b})\},\ i_h\in S_o\\
-\lim_{n\rightarrow \infty}\frac{1}{n}\log P\{\hat{f}(x_{i_h})< \hat{f}(x_{m_b})\},\ i_h\in S_n
\end{aligned}\right\}\\
&=R_{m_b,i_h}(\theta_b,\theta_h, \boldsymbol{\alpha_b},\boldsymbol{\alpha_h})
=\inf_v(\theta_bI_{m_b}(v)+\theta_hI_{i_h}(v))
,\ h\!\neq\! b.\notag
\end{split}
\end{equation}

Since $(\hat{f}(x_{m_b}),\hat{f}(x_{i_h}))$ are all normally distributed random variables, the rate function can be expressed as follows according to the results in \cite{glynn2004large}. Then, we can get
\begin{equation}
\begin{split}
&R_{m_b,i_h}(\theta_b,\theta_h, \boldsymbol{\alpha_b},\boldsymbol{\alpha_h})\\
&=\frac{({f}(x_{m_b})\!-\!{f}(x_{i_h}))^2/2}{\frac{\sigma_{b}^2}{\theta_b}\!\Big(\!\frac{\eta_{m_b\!,1}^2}{\alpha_{1_b}}\!+\!\frac{\eta_{m_b\!,s}^2}{\alpha_{s_b}}\!+\!\frac{\eta_{m_b\!,k}^2}{\alpha_{k_b}}\!\Big)\!+\!\frac{\sigma_{h}^2}{\theta_h}\!\Big(\!\frac{\eta_{i_h\!,1}^2}{\alpha_{1_h}}\!+\!\frac{\eta_{i_h\!,s}^2}{\alpha_{s_h}}\!+\!\frac{\eta_{i_h\!,k}^2}{\alpha_{k_h}}\!\Big)}\!, h\!\neq\! b,
\end{split}\notag
\end{equation}

The proof of (\ref{eq:11a}) is similar to that is given for Lemma 2, and hence is omitted for brevity.

\section*{Proof of Lemma 5}

We first rewrite optimization model (\ref{eq:66}) as
\begin{align}\label{eq.99}
\begin{split}
\max&\ z\\
\text{s.t.}\ &\widetilde{R}_{i_h}(\theta_h, \boldsymbol{\alpha_h})-z\geq0, h=1,...,m, h\neq b, i_h=1_h,...,k_h\\
&R_{i_b}(\theta_b,\boldsymbol{\alpha_b})-z\geq0, i_b=1_b,...,k_{b}\\
&\alpha_{1_h}+\alpha_{s_h}+\alpha_{k_h}=1,h =1,...,m\\
&\sum_{h=1}^l\theta_h=1,\ \theta_h,\alpha_{1_h},\alpha_{s_h},\alpha_{k_h}\geq0, h =1,...,m.
\end{split}
\end{align}
Similarly, the optimization models (\ref{24a}) and (\ref{23a}) can be rewritten as follows:
\begin{align}\label{eq.100}
\begin{split}
\!\!\!\!\max&\ z\\
\text{s.t.}\ &\widetilde{R}_{i_h}(\theta_h, \boldsymbol{\alpha_h})\!-\!z\geq0, i_h\!=\!1_h,...,k_h\\
&\alpha_{1_h}+\alpha_{s_h}+\alpha_{k_h}=1, \alpha_{1_h},\alpha_{s_h},\alpha_{k_h}\geq0,
\end{split}
\end{align}
for $h =1,2,...,m,h\neq b$, and
\begin{align}\label{eq.110}
\begin{split}
\max&\ z\\
\text{s.t.}\
&R_{i_b}(\theta_b,\boldsymbol{\alpha_b})-z\geq0, i_b=1_b,...,k_{b}\\
&\alpha_{1_b}+\alpha_{s_b}+\alpha_{k_b}=1,\alpha_{1_b},\alpha_{s_b},\alpha_{k_b}\geq0.
\end{split}
\end{align}

Let $\boldsymbol{\widetilde{\alpha}_h^*}$ and $\boldsymbol{\widetilde{\alpha}_b^*}$ be the optimal solutions to (\ref{eq.100}) and (\ref{eq.110}).
The optimization model (\ref{eq.99}), (\ref{eq.100}) and (\ref{eq.110}) have the same objective function, and the domain of (\ref{eq.99}) is a subset of (\ref{eq.100}) and (\ref{eq.110}).
Therefore, if $\boldsymbol{\widetilde{\alpha}_h^*}$ and $\boldsymbol{\widetilde{\alpha}_b^*}$ are
feasible to (\ref{eq.99}), they are the optimal solutions of (\ref{eq.99}). According to the formula of $R_{i_b}(\theta_b, \boldsymbol{\alpha_b})$ and $\widetilde{R}_{i_h}(\theta_h, \boldsymbol{\alpha_h})$, we can conclude that $\boldsymbol{\widetilde{\alpha}_h^*}$ and $\boldsymbol{\widetilde{\alpha}_b^*}$ do not depend on the value of $\theta_h$. Therefore, we can conclude that $\boldsymbol{\widetilde{\alpha}_h^*}$ and $\boldsymbol{\widetilde{\alpha}_b^*}$ are feasible to (\ref{eq.99}), and they are optimal to (\ref{eq:66}). Then we get the conclusions in Lemma 5.

\section*{Proof of Theorem 3}

In order to determine the optimal $\boldsymbol{\alpha_h^*}$ for each partition, we consider two cases ($h\neq b$ and $h = b$).

For $h\neq b$, we solve the optimization model (\ref{24a}). Given the definition of the key design $i_h^*$, we rewrite (\ref{24a}) as
\begin{align}\label{eq:10}
\begin{split}
&\max\ \widetilde{R}_{m_b,i_h^*}(\theta_h, \boldsymbol{\alpha_h})\\
\text{s.t.}\
&\alpha_{1_h}+\alpha_{s_h}+\alpha_{k_h}=1,
\alpha_{1_h},\alpha_{s_h},\alpha_{k_h}>0,
\end{split}
\end{align}
for $h=1,2,...,m$ and $h\neq b$.

Since $\widetilde{R}_{i_h^*}(\theta_h, \boldsymbol{\alpha_h})$ is concave and strictly increasing functions of $\boldsymbol{\alpha_h}$ \citep{xiao2015optimal,glynn2004large},
the optimization problem is a convex programming problem. 
We define a Lagrangian function $L\!=\!-\widetilde{R}_{m_b,i_h^*}(\theta_h, \boldsymbol{\alpha_h})+\lambda_h(\alpha_{1_h}\!+\!\alpha_{s_h}\!+\!\alpha_{k_h}-1)$ and determine the optimal allocation based on its Karush-Kuhn-Tucker (KKT) conditions.

\begin{align}
\begin{split}
\widetilde{R}_{m_b,i_h^*}(\theta_h,\boldsymbol{\alpha_h})
\!=\!\frac{(f(x_{m_b})-f(x_{i_h^*}))^2/2}{\frac{\sigma_{h}^2}{\theta_h}\Big(\frac{\eta_{i_h^*,1}^2}{\alpha_{1_h}}+\frac{\eta_{i_h^*,s}^2}{\alpha_{s_h}}+\frac{\eta_{i_h^*,k}^2}{\alpha_{k_h}}\Big)}.
\notag
\end{split}
\end{align}
\begin{align}
\begin{split}
\frac{\partial L}{\partial\alpha_{r_h}}\!=\!-\frac{(f(x_{m_b})-f(x_{i_h^*}))^2/2}
{\frac{\sigma_{h}^2}{\theta_h}\Big(\frac{\eta_{i_h^*,1}^2}{\alpha_{1_h}}+\frac{\eta_{i_h^*,s}^2}{\alpha_{s_h}}+\frac{\eta_{i_h^*,k}^2}{\alpha_{k_h}}\Big)^2}\frac{\eta_{i_h^*,r}^2}{{\alpha_{r_h}^*}^2}\!+\!\lambda_h\!=\!0,\\
{r_h}=1_h,s_h,k_h.
\notag
\end{split}
\end{align}
Using the fact that $\alpha_{1_h}^*+\alpha_{s_h}^*+\alpha_{k_h}^*=1$, we can establish
\begin{equation}\label{24b}
\alpha_{r_h}^*=\frac{|\eta_{i_h^*,r}|}{|\eta_{i_h^*,1}|+|\eta_{i_h^*,s}|+|\eta_{i_h^*,k}|}, {r_h}=1_h,s_h,k_h.
\end{equation}

For $h=b$,
we solve the optimization model (\ref{23a}), which could be rewritten as
\begin{align}\label{eq:23}
\begin{split}
&\max\ {R}_{m_b,i_b^*}(\theta_b, \boldsymbol{\alpha_b})\\
\text{s.t.}\
&\alpha_{1_b}+\alpha_{s_b}+\alpha_{k_b}=1,
\alpha_{1_b},\alpha_{s_b},\alpha_{k_b}>0.
\end{split}
\end{align}
The process of solving the (\ref{eq:23}) is similar to that is given for Theorem 1, and hence is omitted for brevity.
Then, we can get the conclusions in Theorem 3.

\section*{Proof of Theorem 4}

In order to determine the location of the support design $s_h$ for each partition, we consider two cases ($h\neq b$ and $h= b$).

For $h\neq b$, we plug (\ref{eq:9}) into $\widetilde{R}_{m_b,i_h^*}(\theta_h,\boldsymbol{\alpha_h})$, we have
$$\!\widetilde{R}_{m_b,i_h^*}(\theta_h,\boldsymbol{\alpha_h})=\frac{(f(x_{m_b})-f(x_{i_h^*}))^2/2}
{\frac{\sigma_{h}^2}{\theta_h}\left(|\eta_{i_h^*,1}|+|\eta_{i_h^*,s}|+|\eta_{i_h^*,k}|\right)^2}.$$
\begin{align}
\begin{split}
\frac{\partial \widetilde{R}_{m_b,i_h^*}(\theta_h,\boldsymbol{\alpha_h})}{\partial x_{s_h}}&=\frac{-(f(x_{m_b})-f(x_{i_h^*})^2}
{\frac{\sigma_{0h}^2}{\theta_h}\Big(|\eta_{i_h^*,1}|+|\eta_{i_h^*,s}|+|\eta_{i_h^*,k}|\Big)^3}\\
&\times\Big(\frac{\partial|\eta_{i_h^*,1}|}{\partial x_{s_h}}\!+\!\frac{\partial|\eta_{i_h^*,s}|}{\partial x_{s_h}}\!+\!\frac{\partial|\eta_{i_h^*,k}|}{\partial x_{s_h}}\Big).
\notag
\end{split}
\end{align}
where
$$\frac{\partial \eta_{i_h^*,1}}{\partial x_{s_h}}=\frac{(x_{k_h}-x_{i_h^*})(x_{1_h}-x_{i_h^*})}{(x_{1_h}-x_{s_h})^2(x_{1_h}-x_{k_h})},$$
\begin{align}
\begin{split}
\frac{\partial \eta_{i_h^*,s}}{\partial x_{s_h}}
=-\frac{(x_{1_h}-x_{i_h^*})(x_{k_h}-x_{i_h^*})(2x_{s_h}-x_{1_h}-x_{k_h})}{(x_{s_h}-x_{1_h})^2(x_{s_h}-x_{k_h})^2},
\notag
\end{split}
\end{align}
$$\frac{\partial \eta_{i_h^*,k}}{\partial x_{s_h}}=\frac{(x_{1_h}-x_{i_h^*})(x_{k_h}-x_{i_h^*})}{(x_{k_h}-x_{1_h})(x_{k_h}-x_{s_h})^2}.$$
When $x_{s_h}>x_{i_h^*}$,
$$\frac{\partial \widetilde{R}_{m_b,i_h^*}(\theta_h,\boldsymbol{\alpha_h})}{\partial x_{s_h}}<0.$$
When $x_{s_h}<x_{i_h^*}$,
$$\frac{\partial \widetilde{R}_{m_b,i_h^*}(\theta_h,\boldsymbol{\alpha_h})}{\partial x_{s_h}}>0.$$
Therefore, $\widetilde{R}_{m_b,i_h^*}(\theta_h,\boldsymbol{\alpha_h})$ is maximized when
$$x_{s_h}=x_{i_h^*}.$$
According to equation (\ref{eq:9}), we can get
\begin{align}\alpha_{i_h}^*=
\begin{cases}
1,
&i_h=i_h^*,\\
0,
&\text{otherwise}.\notag
\end{cases}
\end{align}
When $i_h^*$ is at the extreme locations ($1_h$ and $k_h$),
\begin{align}
\begin{split}
\widetilde{R}_{m_b,i_h^*}(\theta_h,\boldsymbol{\alpha_h})
\!=\!\frac{(f(x_{m_b})-f(x_{i_h^*}))^2/2}{{\sigma_{h}^2\eta_{i_h^*,r}^2}/{\theta_h}}{\alpha_{r_h}},
r_h=i_h^*.\notag
\end{split}
\end{align}
In this setting, the location of the support design $s_h$ has no effect on $\widetilde{R}_{m_b,i_h^*}(\theta_h,\boldsymbol{\alpha_h})$. For simplicity, we use the D-optimal support design and let $x_{s_h}={(x_{1_h}+x_{k_{h}})/2}$ (round to the nearest design as needed) \citep{kiefer1959optimum}.

For $h=b$, the proof is similar to that is given for Theorem 2, and hence is omitted for brevity.
Then, we can get the conclusions in Theorem 4.

\section*{Proof of Theorem 5}

Given the definition of key design $i_h^*$, we rewrite the optimization model (\ref{eq:6}) as
\begin{align}
\begin{split}
\max\min\Big\{&R_{m_b,i_b^*}(\theta_b,\boldsymbol{\alpha_b^*}),
\min_{1\leq h \leq m\atop h\neq b}\Big(
R_{m_b,i_h^*}(\theta_b, \theta_h,\boldsymbol{\alpha_b^*},\boldsymbol{\alpha_h^*})\Big)\Big\},\\
&\text{s.t.}\ \sum_{h=1}^l\theta_h=1,\ \theta_h\geq0,\ h =1,...,l.
\end{split}\notag
\end{align}

The $R_{m_b,i_b^*}(\theta_b,\boldsymbol{\alpha_b^*})$ and $R_{m_b,i_h^*}(\theta_b, \theta_h,\boldsymbol{\alpha_b^*},\boldsymbol{\alpha_h^*})$ are concave and strictly increasing functions of $\theta_h$ and $\theta_b$ \citep{glynn2004large}. Therefore, the optimization problem is a convex programming problem, and the first order condition is also the optimality condition.

Rewrite the optimization model above as
\begin{align}
\begin{split}
\max&\ z\\
\text{s.t.}\ &R_{m_b,i_h^*}(\theta_b, \theta_h,\boldsymbol{\alpha_b^*},\boldsymbol{\alpha_h^*})-z\geq0, h=1,...,l, h\neq b\\
&R_{m_b,i_b^*}(\theta_b,\boldsymbol{\alpha_b^*})-z\geq0\\
&\sum_{h=1}^l\theta_h=1,\ \theta_h\geq0,\ h =1,...,l.
\end{split}
\end{align}
From the KKT conditions, there exist $\lambda_h\geq 0,\ h=1,2,...,l$ and $\mu>0$ such that
\begin{equation}\label{eq:11}
1-\sum_{h=1}^{m}\lambda_h=0
\end{equation}
\begin{equation}\label{eq:12}
\mu-\lambda_h\frac{\partial R_{m_b,i_h^*}(\theta_b^*, \theta_h^*,\boldsymbol{\alpha_b^*},\boldsymbol{\alpha_h^*})}{\partial\theta_h}=0, h=1,2,..,l,h\neq b
\end{equation}
\begin{equation}
\mu-\!\sum_{h=1\atop h\neq b}^l\lambda_h\frac{\partial R_{m_b,i_h^*}(\theta_b^*, \theta_h^*,\boldsymbol{\alpha_b^*},\boldsymbol{\alpha_h^*})}{\partial\theta_b}-\lambda_b\frac{\partial R_{m_b,i_b^*}(\theta_b^*,\boldsymbol{\alpha_b^*})}{\partial\theta_b}=0
\end{equation}
\begin{equation}\label{eq:14}
\lambda_h(z-R_{m_b,i_h^*}(\theta_b^*, \theta_h^*,\boldsymbol{\alpha_b^*},\boldsymbol{\alpha_h^*}))=0
, h=1,2,..,l,h\neq b
\end{equation}
\begin{equation}\label{eq:15}
\lambda_b(z-R_{m_b,i_b^*}(\theta_b^*,\boldsymbol{\alpha_b^*}))=0
\end{equation}

From (\ref{eq:11}), it can be concluded that there must exist some $\lambda_h>0, h=1,2,...,l$. If there is one $\lambda_h=0,h=1,2,...,l, h\neq b$, it results $\mu=0$ based on (\ref{eq:12}). That means all the $\lambda_h=0$, as ${\partial R_{m_b,i_h^*}(\theta_b, \theta_h,\boldsymbol{\alpha_b^*},\boldsymbol{\alpha_h^*})}/{\partial\theta_h}$ is strictly positive. Therefore, we can conclude that $\lambda_h\neq 0,\ h=1,2,...,l, h\neq b.$ According to (\ref{eq:14}),
\begin{equation}\label{eq:30}
z^*=R_{m_b,i_h^*}(\theta_b^*, \theta_h^*,\boldsymbol{\alpha_b^*},\boldsymbol{\alpha_h^*}),\ h=1,2,...,l, h\neq b.
\end{equation}

In addition, since $\theta_h/\theta_b\rightarrow0$ as $l$ goes to infinity and $R_{m_b,i_b^*}(\theta_b^*,\boldsymbol{\alpha_b^*})>R_{m_b,i_h^*}(\theta_b^*, \theta_h^*,\boldsymbol{\alpha_b^*},\boldsymbol{\alpha_h^*}),\ h=1,2,...,l, h\neq b$, it can be concluded that $\lambda_b=0$ according to (\ref{eq:15}). Therefore, we have
\begin{align}\label{eq.40}
\begin{split}
&\sum_{h=1\atop h\neq b}^l\!\frac{\partial\! R_{m_b,i_h^*}\!(\theta_b^*,\! \theta_h^*,\!\boldsymbol{\alpha_b^*},\!\boldsymbol{\alpha_h^*})/\partial\theta_b}{\partial\! R_{m_b,i_h^*}\!(\theta_b^*,\! \theta_h^*,\!\boldsymbol{\alpha_b^*},\!\boldsymbol{\alpha_h^*})/\partial\theta_h}\!=1.
\end{split}
\end{align}
According to Theorem 4, we have
\begin{align}
\begin{split}
R_{m_b,i_h^*}(\theta_b^*,\theta_h^*,\boldsymbol{\alpha_b^*},\boldsymbol{\alpha_h^*})=
\frac{(f(x_{m_b})-f(x_{i_h^*}))^2/2}{\frac{\sigma_{b}^2}{\theta_b^*}\Big(\frac{\eta_{m_b,1}^2}{\alpha^*_{1_b}}+\frac{\eta_{m_b,s}^2}{\alpha^*_{s_b}}+\frac{\eta_{m_b,k}^2}{\alpha^*_{k_b}}\Big)+\frac{\sigma_{h}^2}{\theta_h^*}}.
\end{split}
\end{align}
Then, we can get
\begin{align}
\begin{split}
&\frac{\partial R_{m_b,i_h^*}(\theta_b^*,\theta_h^*,\boldsymbol{\alpha_b^*},\boldsymbol{\alpha_h^*})}{\partial \theta_b^*}\\
&=\frac{\sigma_{b}^2(f(x_{m_b})-f(x_{i_h^*}))^2\Big(\frac{\eta_{m_b,1}^2}{\alpha^*_{1_b}}+\frac{\eta_{m_b,s}^2}{\alpha^*_{s_b}}+\frac{\eta_{m_b,k}^2}{\alpha^*_{k_b}}\Big)/{\theta^*_b}^2}{2\bigg(\frac{\sigma_{b}^2}{\theta_b^*}\Big(\frac{\eta_{m_b,1}^2}{\alpha^*_{1_b}}+\frac{\eta_{m_b,s}^2}{\alpha^*_{s_b}}+\frac{\eta_{m_b,k}^2}{\alpha^*_{k_b}}\Big)+\frac{\sigma_{h}^2}{\theta_h^*}\bigg)^2},
\end{split}
\end{align}
and
\begin{align}
\begin{split}
\frac{\partial R_{b_m,i_h^*}(\theta_b^*,\theta_h^*,\boldsymbol{\alpha_b^*},\boldsymbol{\alpha_h^*})}{\partial \theta_h^*}=
\!\frac{\sigma_{h}^2(f(x_{m_b})-f(x_{i_h^*}))^2/{\theta^*_h}^2}{2\bigg(\frac{\sigma_{b}^2}{\theta_b^*}\Big(\frac{\eta_{m_b,1}^2}{\alpha^*_{1_b}}+\frac{\eta_{m_b,s}^2}{\alpha^*_{s_b}}+\frac{\eta_{m_b,k}^2}{\alpha^*_{k_b}}\Big)+\frac{\sigma_{h}^2}{\theta_h^*}\bigg)^2},
\end{split}
\end{align}


As $l$ goes to infinity, plug the results above into (\ref{eq:30}) and
(\ref{eq.40}), we can get the equations in Theorem 5.

\end{document}